\pdfoutput=1
\ifx \forjournal\undefined
\documentclass[nologo]{uq-report}
\else
\documentclass[nologo,notitlepageheader]{dakota-article}
\fi

\shortauthor{Jakeman, Franzelin, Narayan, Eldred, Plf\"{u}ger}
\usepackage{amsmath,amsfonts,amssymb,graphicx,url,amsthm,hyperref,mathtools,algorithmic,algorithm}
\def\SNL{Optimization and Uncertainty Quantification, Sandia National Laboratories, Albuquerque, NM, 87123}
\def\STUTTGART{Department of Computer Science, University of Stuttgart, Stuttgart, Germany}
\def\UTAH{Department of Mathematics and Scientific Computing and Imaging (SCI) Institute, University of Utah,  Salt Lake City, UT, USA}
\corauthor{J.D. Jakeman}
\coremail{jdjakem@sandia.gov}
\title{Polynomial chaos expansions for dependent random variables}
\shorttitle{PCE for dependent random variables}

\def\rv{Z}  
\def\rvd{z}  
\def\brv{Z} 
\def\pbwt{\omega} 
\newcommand{\norm}[2]{\left\lVert #1 \right\rVert_{#2}}

\def\R{\mathbb{R}}
\def\N{\mathbb{N}}
\def\nodes{\mathcal{Z}}

\newcommand{\rvdsamp}[1]{\rvd^{(#1)}}
\DeclareMathOperator*{\argmax}{argmax}

\newcommand{\dx}[1]{\mathrm{d}#1}

\def\urv{U}
\def\burv{U}
\def\urvd{u}

\def\Rvx{R^{\brv}}
\def\Rvy{{R^{V}}}

\def\dom{\Omega}
\def\rotation{A}

\date{\today}

\ifx \forjournal\undefined
\newcommand{\rev}[1]{{#1}}

\usepackage{authblk}

\author[1]{J.D. Jakeman}
\affil[1]{\SNL} 
\author[2]{F. Franzelin}
\affil[2]{\STUTTGART}
\author[3]{A. Narayan}
\affil[3]{\UTAH}
\author[1]{M.S. Eldred}
\author[2]{Dirk Plf\"{u}ger}
\def\gaussicdf{{\Phi_g}} 
\def\ngs{J} 

\else
\funding{See acknowledgements}
\newtheorem{definition}{Definition}[section]
\author{John D. Jakeman\thanks{\SNL}, Fabian Franzelin\thanks{\STUTTGART}, Akil Narayan\thanks{\UTAH}, Michael Eldred\samethanks[1], Dirk Plf\"{u}ger\samethanks[2]}
\date{\today}
\hypersetup{
  pdfkeywords={},
  pdfsubject={},
  pdfcreator={Emacs 25.2.2 (Org mode 8.2.10)}}

\newcommand{\rev}[1]{{\color{blue}#1}}

\def\gaussicdf{{\color{blue}\Phi_g}} 
\def\ngs{{\color{blue}{J}}} 
\fi

\begin{document}

\maketitle
\begin{abstract}
  Polynomial chaos expansions (PCE) are well-suited to quantifying uncertainty in models parameterized by independent random variables. The assumption of independence leads to simple strategies for building multivariate orthonormal bases and for sampling strategies to evaluate PCE coefficients. In contrast, the application of PCE to models of dependent variables is much more challenging. Three approaches can be used to construct PCE of models of dependent variables. The first approach uses mapping methods where measure transformations, such as the Nataf and Rosenblatt transformation, can be used to map dependent random variables to independent ones; however we show that this can significantly degrade performance since the Jacobian of the map must be approximated. A second strategy is the class of dominating support methods. In these approaches a PCE is built using independent random variables whose distributional support dominates the support of the true dependent joint density; we provide evidence that this approach appears to produce approximations with suboptimal accuracy. A third approach, the novel method proposed here, uses Gram-Schmidt orthogonalization (GSO) to numerically compute orthonormal polynomials for the dependent random variables. This approach has been used successfully when solving differential equations using the intrusive stochastic Galerkin method, and in this paper we use GSO to build PCE using a non-intrusive stochastic collocation method. The stochastic collocation method treats the model as a black box and builds approximations of the input-output map from a set of samples. Building PCE from samples can introduce ill-conditioning which does not plague stochastic Galerkin methods. To mitigate this ill-conditioning we generate weighted Leja sequences, which are nested sample sets, to build accurate polynomial interpolants. We show that our proposed approach, GSO with weighted Leja sequences, produces PCE which are orders of magnitude more accurate than PCE constructed using mapping or dominating support methods.
\end{abstract}

\section{Introduction}
\label{sec-1}
All models of modern scientific applications are subject to various sources of uncertainty. The effect of these uncertainties on model predictions can be assessed by viewing the model as an input-output map, where the inputs are random variables with known (but possibly complicated) distribution, and the output is a set of model quantities of interest (QoI). Uncertainty quantification (UQ) then refers to the process of computing the output statistics that result from this input-output map.

UQ of high-fidelity models typically requires large numbers of simulations. Building an approximation or surrogate for the model input-output map is an effective and popular approach to reduce the computational burden of UQ, and numerous techniques have been developed to build such approximations. Within the computational science and engineering community, some of the most widely adopted methods for approximating models parameterized by random variables are those based on generalized polynomial chaos expansions \cite{Ghanem_book_1991,Xiu_K_SISC_2002}, sparse grid approximation \cite{Xiu_H_SISC_2005,Nobile_TW_SINUM_08}, Gaussian process models \cite{rasmussen2005} and low-rank tensor decompositions \cite{Oseledets2011}. These methods can be very efficient when building approximations of models subject to independent random variables, and as such, most computational tools are built assuming the input variables are independent.  However, there is a dearth of algorithmic options when the variables are dependent.

In this paper we present sampling strategies for building polynomial chaos expansions (PCE) of models influenced by a high-dimensional random vector that are accurate and stable when the components of the random vector are not independent. PCE essentially seeks to build a polynomial approximation of a function (model) whose inputs are the random variables; our goal will be the generation of a PCE surrogate that is accurate in a norm weighted by the probability density function of the dependent variables. PCE are ideally suited to approximating functions of random variables because they employ basis functions which are orthonormal to the probability measure \(\pbwt\) of the variables, and this relationship can be exploited to construct stable approximation schemes that are accurate in regions of high probability.

The stochastic Galerkin \cite{Ghanem_book_1991} and stochastic collocation \cite{Babuska_NT_SNA_2007,Xiu_H_SISC_2005} methods are the two main approaches for generating a PCE surrogate, which amounts to computing a set of PCE coefficients. In this paper we focus on stochastic collocation because it allows the computational model to be treated as a black box. Stochastic collocation proceeds in two steps: (i) running the computational model with a set of realizations of the random parameters and (ii) constructing an approximation of the corresponding model output. In the relevant situation when the model is expensive, step (i) is the most time-consuming portion. Stochastic collocation methods include pseudospectral projection \cite{Conrad_M_SISC_2013,Constantine_EP_CMAME_2012}, sparse grid interpolation \cite{Barthelmann_NR_ACM_2000,Buzzard_RESS_2012}, least orthogonal interpolation \cite{Narayan_X_SISC_2012}, least squares \cite{Migliorati_NSTSISC_2013,Tang_Z_SISC_2014}, compressive sensing \cite{Doostan_O_JCP_2011,Yan_GX_IJUQ_2012,Jakeman_ES_JCP_2015} and low-rank tensor decompositions \cite{Chevreuil_LNR_SIAMJUQ_2015,Doostan_VI_CMAME_2013,Gorodetsky_J_JCP_2018}.  

Sparse grid interpolations and pseudospectral approximations are by construction ideally suited to approximation when the input random variables are independent. Polynomial-based sparse grid methods utilize univariate high-order (e.g., Gaussian) quadrature rules that are optimal for tensor product approximation \cite{Nobile_TW_SINUM_08,Narayan_J_SISC_2014}. Recently, sub-sampled tensor-product quadrature has been used to generate sample sets for independent random variables in the context of interpolation \cite{Li_Z_WRR_2007}, least squares \cite{Seshadri_NM_SIAMUQ_2017} and sparse regression \cite{Tang_I_SIAMUQ_2014}. There has also been extensive work on using random sampling for least squares \cite{Migliorati_NSTSISC_2013,Narayan_JZ_MC_2017} and compressive sensing \cite{Rauhut_W_JAT_2012,Hampton_D_JCP_2015,Jakeman_NZ_SISC_2017} as well as some work on generating deterministic sequences for interpolation \cite{Narayan_X_SISC_2012,Narayan_J_2018}.

The aforementioned methods usually rely on tensor product structure and are most effective for independent random variables. If one is only interested in computing moments, recent advances in polynomial quadrature for dependent measures can be utilized \cite{Arnst_GPR_IJNME_2012,Constantine_PW_IJNME_2014,keshavarzzadeh_numerical_2017,cui_stochastic_2018}.
The polynomial accuracy of these methods is inadequate for efficient pseudo spectral projection. Consequently regression-based methods are the only viable alternative. 

If $f$ is the model under consideration and $Z$ is a vector of random parameters that are input to the model, a PCE approach posits the representation
\begin{equation}
\label{eq:pce-integer-index}
f(z)=\sum_{n=1}^\infty\alpha_{n}\phi_{n}(z),
\end{equation}
where $z$ are realizations of the variables $Z$, $\alpha_n$ are the PCE coefficients that must be computed, and the basis functions $\phi_n$ are polynomial basis functions that are pairwise orthonormal under an inner product weighted by the probability density of $Z$. Above we assume that $f$ is scalar-valued, but the procedures we describe carry over to vector- or function-valued outputs.
When the components of $Z$ are independent, one can generate the multivariate polynomials $\phi_n$ from univariate orthogonal polynomials, but such a construction is not easily accomplished when $Z$ has dependent components. One approach for dealing with dependent variables is to build an approximation for a set of independent variables whose tensor-product measure dominates the dependent measure~\cite{Jakeman_EX_JCP_2010,Chen_PX_JCP_2013}. Such an approach introduces an error which is proportional to the ``distance'' between the tensor-product and dependent measures.

In this paper we use Gram-Schmidt orthogonalization to generate a set of polynomials orthonormal to arbitrary probability measures. This approach was first proposed in the multivariate setting in \cite{Witteveen_B_AIAA_2006} for stochastic Galerkin methods and used for solving time-dependent PCE using stochastic Galerkin projection \cite{Gerritsma_SVK_2010}. However, the use of Gram-Schmidt orthogonalization has received comparatively less attention for stochastic collocation. We will observe in this paper a well-known phenomenon, that the use of a Gram-Schmidt procedure to build basis functions is poorly conditioned. However we propose a method that can reduce the amount of ill-conditioning by using specialized preconditioning and sampling approaches. Specifically, in this paper we adapt and improve the strategies for building weighted Leja sequences for interpolation developed for independent random variables in \cite{Narayan_J_2018}.


This paper is devoted to computational studies of novel and recently-developed algorithms for computing PCE expansions for dependent variables. The novelty and main outcomes of this study are described below:
\begin{itemize}
  \item We provide a systematic study, on practical problems, which compares the use of mapping methods (such as the Rosenblatt transformation), measure domination methods, and Gram-Schmidt Orthogonalization (GSO) methods. In short, we observe that GSO provides the most accurate procedure, but may be ill-conditioned in some cases. However, this ill-conditioning frequently does not degrade the accuracy of the resulting emulator.
  \item We propose the use of weighted Leja sequences (where the weight is related to the dependent density) when constructing PCE approximations for dependent measures. An ingredient in this approach is the use of GSO for constructing an orthonormal basis for the dependent measure. We observe that this approach frequently performs much better than the alternatives described above. 
\end{itemize}

In the remainder of this paper we introduce polynomial chaos expansions and their construction for both tensor-product and dependent probability measures. We then propose stable and accurate sampling schemes for regression and interpolation, using PCE, which can be used for arbitrary measures. We then provide an investigation of the Gram-Schmidt procedure, compare its performance with probabilistic transformation maps for multivariate approximation, and conclude with a number of numerical examples to highlight the strengths of our proposed approach.

\section{Polynomial chaos expansions for independent random variables}
\label{sec-2}
Let \(f: \R^d \rightarrow \R\) be a function (model) of a \(d\)-variate random variable \(\brv=(\rv_1,\ldots,\rv_d)\). The random variable has associated probability density function \(\pbwt(\rvd)\) for \(\rvd\in\dom \subset \R^d\).  Polynomial chaos expansions represent the model output \(f(\rvd)\) as an expansion in orthonormal polynomials, as in \eqref{eq:pce-integer-index}.  The basis functions $\phi_n$ are typically constructed to be orthonormal with respect to the density \(\pbwt\), that is 
\[
  (\phi_{i}(\rvd),\phi_{j}(\rvd))_{L^2_{\pbwt}(\Omega)} \coloneqq \int_{\dom} \phi_{i}(\rvd)\phi_{j}(\rvd) \pbwt(\rvd)= \delta_{i,j},
\]
where \(\dom\) is the support of the density $\pbwt$ and $\delta_{i,j}$ is the Kronecker delta function. Under mild conditions on the distribution \(\pbwt(\rvd)\), any function $f(\rvd)$ with finite variance, i.e. \(f\in L^2_{\pbwt}(\Omega)\) can be represented by a PCE that converges in \(L^2_{\pbwt}\) to the true function asymptotically \cite{ernst_convergence_2012}.  In this section we describe the ``canonical" construction when $Z$ has independent components, which is well-known and standard. The next section discusses how one can construct a PCE basis when $Z$ has dependent components. 

Polynomial chaos expansions are most easily constructed when the components of \(\rv\) are independent.  Under the assumption of independence, we have
\begin{align*}
  \dom &= \times_{i=1}^d \dom_i, & \dom_i &\subset \R, & \pbwt(\rvd) &= \prod_{i=1}^d \pbwt_i(\rvd_i),
\end{align*}
where $\pbwt_i$ are the marginal densities of the variables \(\brv_i\), which completely characterizes the distribution of $\rv$. This allows us to express the basis functions \(\phi\) as tensor products of univariate orthonormal polynomials. That is
\begin{align}
\phi_\lambda(\rvd)=\prod_{i=1}^d \phi^i_{\lambda_i}(\rvd_i),
\end{align}
where \(\lambda=(\lambda_1\ldots,\lambda_d)\in\mathbb{N}_0^d\) is a multi-index, and the univariate basis functions $\phi^i_j$ are defined uniquely (up to a sign) for each $i = 1, \ldots, d$, as
\begin{align*}
  \int_{\Omega_i} \phi^i_{j}(z_i) \phi^i_{k}(z_i) \pbwt_i(z_i) \dx{z_i} &= \delta_{j,k}, & j, k &\geq 0, & \deg \phi^i_j &= j.
\end{align*}
In practice the PCE \eqref{eq:pce-integer-index} must be truncated to some finite number of terms, say \(N\), defined by a multi-index set $\Lambda \subset \N_0^d$:
\begin{align}
\label{eq:pce-multi-index}
f(\rvd) &\approx f_N(\rvd) = \sum_{\lambda\in\Lambda}\alpha_{\lambda}\phi_{\lambda}(\rvd), & |\Lambda| &= N.
\end{align}
Frequently the PCE is truncated to retain only the multivariate polynomials whose associated multi-indices have norm at most \(p\), i.e.,
\begin{align}
\label{eq:hyperbolic-index-set} 
\Lambda &= \Lambda^d_{p,q} = \{\lambda \mid \norm{\lambda}{q} \le p\}., & \left\| \lambda \right\|_q &\coloneqq \left(\sum_{i=1}^d \lambda^q_i\right)^{1/q}.
\end{align}
Taking \(q=1\) results in a total-degree space having dimension \(\text{card}\; \Lambda^d_{p,1} \equiv N = { d+p \choose d }\). The choice of $\Lambda$ identifies a subspace in which $f_N$ has membership:
\begin{align*}
  \pi_\Lambda &\coloneqq \mathrm{span} \left\{ \phi_\lambda \;\; \big| \;\; \lambda \in \Lambda \right\}, & f_N &\in \pi_\Lambda.
\end{align*}
Under an appropriate ordering of multi-indices, the expression \eqref{eq:pce-multi-index}, and the expression \eqref{eq:pce-integer-index} truncated to the first $N$ terms, are identical. Defining \([N]:=\{1,\ldots,N\}\), for \(N\in\mathbb{N}\), we will in the following frequently make use of a linear ordering of the PCE basis, \(\phi_k\) for \(k \in [N]\) from \eqref{eq:pce-integer-index}, instead of the multi-index ordering of the PCE basis \(\phi_{\lambda}\) for \(\lambda \in \Lambda\) from \eqref{eq:pce-multi-index}.  Therefore, 
\begin{align*}
  \sum_{\lambda \in \Lambda} \alpha_\lambda \phi_\lambda(z) = \sum_{n=1}^N \alpha_n \phi_n(z).
\end{align*}
Any bijective map between $\Lambda$ and $[N]$ will serve to define this linear ordering, and the particular choice of this map is not relevant in our discussion.

\section{PCE for dependent random variables}
\label{sec-2-1}
In this section we present three of the most popular strategies for constructing a PCE surrogate when the distribution of $Z$ is not of tensor-product form. The first strategy is a mapping method, where the dependent coordinates are mapped to a set of independent coordinates. The second algorithm identifies a dominating measure of tensor-product form, and a PCE is built with respect to the dominating measure. The final strategy attempts to explicitly compute an orthonormal polynomial basis in $L^2_\omega$ and to subsequently construct a PCE approximation. 

Much of our notation from the previous section carries over, except that we require more caution when speaking of multi-index sets $\Lambda$. Without a tensor-product construction, the definition of the polynomial $\phi_\lambda$ is not unique, and thus the meaning of the subspace $\pi_\Lambda$ is unclear. However, the following definition mitigates part of this issue:
\begin{definition}
  Let $\Lambda \subset \N_0^d$ be a downward-closed index set, i.e., $\lambda \in \Lambda$ implies that $\nu \in \Lambda$ for all $\nu \in \N_0^d$ satisfying $\nu \leq \lambda$. Here, $\leq$ is the partial lexicographic ordering on the $d$-dimensional integer lattice. Then we define
  \begin{align*}
    \pi_\Lambda \coloneqq \mathrm{span} \left\{ \prod_{i=1}^d z_i^{\lambda_i} \;\; \big| \;\; \lambda \in \Lambda \right\}.
  \end{align*}
\end{definition}
This definition is consistent with the meaning of $\pi_\Lambda$ for independent variables introduced in Section~\ref{sec-2} (assuming $\Lambda$ is downward-closed).

\subsection{Mapping methods}
\label{sec-4}
Let \(\mathcal{T}: \dom \rightarrow \rev{\tilde{\dom}}\) denote an invertible transformation which maps possibly dependent random variables \(\brv\) to independent random variables \(\mathcal{T}(Z) = \burv=(\urv_1,\ldots,\urv_d)\) with marginal densities \(\rho_i(u_i)\), \(i\in[d]\). Typically, mapping methods are defined so that $U$ has a uniform distribution on $\rev{\tilde{\dom}}=[0,1]^d$ (and hence is a tensor-product distribution). In this case, we can form the approximation $f_N$ by first forming an approximation in $U$ space~\cite{Eldred_B_AIAA_2009,Torre_MES_ARXIV_2017}:
\begin{align}
\label{eq:pce-u-space}
f(\mathcal{T}^{-1}(\urvd)) \eqqcolon g(\urvd) &\approx g_N(\urvd) = \sum_{\lambda \in \Lambda} \alpha_{\lambda} \psi_{\lambda}(\urvd), & f_N(\rvd) &\coloneqq g_N(\mathcal{T}(\rvd))
\end{align}
Here, \(\psi\) is a tensor product basis of univariate polynomials orthonormal with respect to the univariate densities \(\rho_i\), in this case the univariate orthonormal Legendre polynomials on $[0,1]$. 

This approach is straightforward to implement if $\mathcal{T}$ is available, but the construction of $\mathcal{T}$ is often the bottleneck in implementations. In the following we describe two popular approaches for computing the transformation $\mathcal{T}$.

\subsubsection{Rosenblatt Transformation}
\label{sec-4-1}
When the joint cumulative distribution function  \(F_\brv\) of $Z$ is continuous, an explicit construction of \(\mathcal{T}\) is given by the Rosenblatt transformation \cite{Rosenblatt_AMS_52}. The Rosenblatt transformation defines the following components for the transformation $\mathcal{T}$: 
\begin{align}
  \urvd_1&=F_1(\rvd_1), &\urvd_2&=F_{2\mid 1}(\rvd_2\mid \rv_1 = \rvd_1)  &\cdots & &\urvd_d&=F_{d|d-1, \dots, 1}(\rvd_d\mid \rv_1 = \rvd_1,\ldots,\rv_{d-1} = \rvd_{d-1}),
\end{align}
where \(F_{i|i-1, \dots, 1}(\cdot \mid \rv_1 = \rvd_1,\ldots,\rv_{i-1} = \rvd_{i-1})\) is the conditional distribution of $\rv_i$ given that $(\rv_1, \ldots, \rv_{i-1}) = (\rvd_1, \ldots, \rvd_{i-1})$. 
The marginalization needed to compute the conditional distributions generally requires multivariate integration. For example, computing \(\urvd_1\) requires integration over \(d-1\) variables, and hence a quadrature approach suffers the curse of dimensionality.

The inverse Rosenblatt transformation can be used to obtain the original sample \(\brv\) from the  sample \(\burv\) by solving the following optimization problem
\begin{equation}
  \label{eq:inverse-rosenblatt-transformation}
  \begin{aligned}[t]
    \urv_1 - F_1(\rvd_1) &= 0, &
    \urv_2 - F_{2|1}(\rvd_2|\rv_1=\rvd_1) &= 0, &
    \ldots & &
    \urv_d - F_{d|d-1, \dots, 1}(\rvd_d|\rv_{d-1}=\rvd_{d - 1}, \dots, \rv_{1}=\rvd_1) &= 0.
  \end{aligned}
\end{equation}
For a fixed order of the one-dimensional transformations, there exists
a unique solution to this problem, since the \(F_i\) are strictly monotonic increasing functions. We can solve this problem by applying standard root-finding algorithms, such as the bisection method.

The Rosenblatt transformation essentially requires a closed form for the joint density \(\pbwt(\brv)\). This requirement can be relaxed by utilizing density estimation techniques, such as Gaussian kernel or sparse grid density estimation to construct an approximation to the joint density \cite{Franzelin_Thesis_2018}, which can help mitigate the curse of dimensionality. 

\subsubsection{Nataf Transformation}
\label{sec-4-2}
When the dependence between the random variables \(\brv\) is linear then the Rosenblatt transformation can be simplified. The resulting Nataf transformation~\cite{Liu_K_PEM_1986} requires the availability of the correlation matrix \(\Rvx \in \R^{d \times d}\), of $Z$, the marginal distributions \(F_i(\rvd_i)\), and the marginal density of $\rv_i$, denoted $\pbwt_i(\rvd_i)$.

The Nataf transformation assumes that the joint density \(\pbwt(\brv)\) can be expressed in terms of the multivariate Gaussian density
\begin{equation}
  \label{eq:multivariate-normal-density}
  \eta_{R^V}(\urvd) = \frac{1}{\sqrt{(2 \pi)^d \mathrm{det} (\Rvy)}}
    \exp\left(-\frac{1}{2}\urvd^T \Rvy^{-1} \urvd\right)\;,
\end{equation}
with correlation matrix $R^V$ that must be computed. The non-tensorial density $\pbwt$ is assumed to have the form
\begin{equation}
\label{eq:nataf-density}
  \pbwt(\rvd) = \frac{\eta_{R^V}(\hat{u})}{\prod_{i=1}^d\eta(\hat{u}_i)} \prod_{i=1}^d \pbwt_i(\rvd_i)\;.
\end{equation}
Here \(\pbwt_i\) are the marginal distributions of the joint density \(\pbwt\) and \(\hat{u}_i=\gaussicdf^{-1}(F_i(\rvd_i))\), with \(\eta\) and \(\gaussicdf\) respectively denoting the univariate standard normal density and cumulative distribution functions. Each entry \(\Rvx_{ij}\) of the correlation matrix \(\Rvx\) is related to \(R^V_{ij}\) by
\begin{equation}
  \label{eq:correlation-coefficient-connected}
  \Rvx_{ij} = \int_\R \int_\R \left(\frac{\rvd_i - \mathbb{E}(\rv_i)}{\sqrt{\text{Var}(\rv_i)}}\right)
  \left(\frac{\rvd_j - \mathbb{E}(\rv_j)}{\sqrt{\text{Var}(\rv_j)}}\right) \frac{\eta_{R^{V_i,V_j}}(\gaussicdf^{-1}(F_i(\rvd_i)), \gaussicdf^{-1}(F_j(\rvd_j)))}{\eta(\gaussicdf^{-1}(F_i(\rvd_i))) \eta(\gaussicdf^{-1}(F_j(\rvd_j)))} \pbwt_i(\rvd_i)\pbwt_j(\rvd_j)\dx{\rvd_i}\dx{\rvd_j},
\end{equation}
where
\(R^{V_i,V_j} =\begin{pmatrix}
    1 & R^V_{ij} \\
    R^V_{ij} & 1\;
\end{pmatrix} \in \R^{2 \times 2}.\)
The integral in \eqref{eq:correlation-coefficient-connected} is two-dimensional and so can be solved
with Gaussian quadrature and standard root finding algorithms, such as bisection
\cite{Li_LY_CSB_2008}.

Applying the Nataf transformation consists of three steps. First the ``corrected'' correlation matrix \(\Rvy\) must be constructed by solving \eqref{eq:correlation-coefficient-connected}. This matrix is independent of \rev{the marginals of} the variable \(\brv\) and so can be computed once and stored for any subsequent use. The second step creates a set of intermediate correlated Gaussian random variables \(\hat{U}=(\hat{U}_1,\ldots,\hat{U}_d)\) created from \(\brv\) via \(\hat{U}_i = \gaussicdf^{-1}(F_i(\rv_i))\). Finally these intermediate variables are decorrelated using \(\burv = L^{-1}\hat{U}\), where \(L\) is the lower triangular Cholesky factor of \(\Rvy\). In this paper we \rev{use this Nataf trafansformation \(\mathcal{T}_\text{nataf}^\text{gauss}\) to build PCE using univariate Hermite polynomials. We will also use the transformation $\label{eq:nataf+}\mathcal{T}_\text{nataf}^\text{unif}=\mathcal{T}_{n2u}\circ\mathcal{T}_\text{nataf}^\text{gauss},$
where \(\mathcal{T}_{n2u}\) maps independent normally-distributed random variables to independent uniform variables on \([-1,1]\). When using this transformation, we will build PCE based upon tensor-products of univariate Legendre polynomials.}

There is a close connection between the Nataf transformation and copula dependence modeling. Sklar's theorem states that any joint distribution can be expressed as the product of univariate marginal distributions and a copula that describes the dependency structure. The Nataf density in \eqref{eq:nataf-density} does exactly this, using the Gaussian copula \(\eta_{R^V}(V)/\prod_{i=1}^d\eta_i(V_i)\). We will utilize this connection between the Nataf transformation and Gaussian copula modeling to efficiently generate dependent multivariate samples to test the performance of our proposed method in Section \ref{sec-5}.  However we will show in the same section that the usage of Nataf transformations for constructing approximations via \eqref{eq:pce-u-space} produces poor approximations when compared to the other approaches considered in this paper.

\subsection{Domination methods}\label{ssec:pce-domination}
A domination method can be used to compute a PCE when $Z$ has dependent components.  This approach uses PCE's consisting of tensor product basis functions to approximate functions of dependent random variables, but this approximation comes at a cost \cite{Xiu_K_SISC_2002,Jakeman_EX_JCP_2010,Chen_PX_JCP_2013}. Given $\dom$ and $\pbwt$, the basic idea is to identify a tensor-product density $g$ with support $G \subset \R^d$ so that $G \supseteq \Omega$. For example, if $\Omega$ is compact then one may identify $G$ as the smallest bounding box for $\dom$, and define $g$ as the uniform probability density over $G$. One then constructs a PCE using the density $g$ with the strategy from Section \ref{sec-2}. The error committed by this strategy is essentially well-understood.

Denote the \(g\)-weighted \(L^p\) norm as
\begin{align}
  \lVert f \rVert_{L^p_g} &= \left(\int_{G}\lvert f\rvert^p g(x)\,dx\right)^{1/p},& 1&\le p<\infty.
\end{align}
We have the following lemma, which characterizes accuracy in a $\pbwt$-weighted norm given an approximation that is accurate in the $g$-weighted norm.
\begin{lemma}[Strong convergence \cite{Chen_PX_JCP_2013}]
\label{lem:chen-2013-strong}
Let $p$ and $q$ be conjugate exponents, i.e., \(1/p+1/q=1\) with $p, q \geq 1$.  Assume the error of the approximation \(f_N\) of \(f\) satisfies 
\begin{align}
\label{eq:lp-error-assumption}
\epsilon_N &= \lVert f-f_N \rVert_{L_g^p}, & p &\ge 1.
\end{align}
Then,
\begin{align}
\label{eq:episetmic-theorem-result}
\lVert f-f_N\rVert_{L^p_{\pbwt}} &\le C_{r}^{1/p}\epsilon_N, & C_r \coloneqq \max_{z \in \Omega} \frac{\pbwt(z)}{g(z)},
\end{align}
provided $C_r < \infty$.
\end{lemma}
This lemma states that the error induced by approximating using a polynomial construction that is accurate in a dominating measure induces an error that grows as the ``difference'' between the optimal and non-optimal basis increases. This difference is quantified as the maximum $C_r$ of the ratio of the two densities over $\dom$. Thus if we use a PCE based upon the tensor product of univariate orthogonal polynomials for approximating a function of highly dependent variables then this can induce a substantial increase in error for a fixed sample size, when compared to the error obtained using a polynomial approximation that is designed to be accurate in the original $\pbwt$-weighted norm. We demonstrate in Section \ref{sec-5} that this happens in practice.

\subsection{Orthogonalization methods}\label{ssec:pce-orthogonalization}
The most straightforward way to construct a PCE basis for dependent variables uses Gram-Schmidt orthogonalization \cite{Witteveen_B_AIAA_2006}. Assume any linearly independent set of polynomials \(\{\psi_n\}_{n=1}^N\) is given (a common choice might be polynomials that are orthonormal under a related tensorial measure). Then the polynomials $\phi_n$ orthonormal under $\pbwt$ are constructed numerically by setting \(\phi_1(\rvd)\equiv 1\) and computing
\begin{align}\label{eq:gram-schmidt-alg}
  \phi_n(\brv) &= \frac{\widetilde{\phi}_n}{\|\widetilde{\phi}_n\|_{L^2_\pbwt}}, & \widetilde{\phi}_n(\rvd) &= \psi_n(\rvd)-\sum_{k=1}^{n-1}\frac{(\psi_n,\phi_k)_{L^2_{\pbwt}}}{(\phi_k,\phi_k)_{L^2_{\pbwt}}}\phi_k(\rvd), \quad n\in[N]. 
\end{align} 
The $L^2_\omega$ norms and inner products above can be computed approximately using quadrature. Let \(\nodes_\ngs=(\rvdsamp{1},\ldots,\rvdsamp{\ngs})\), with \(\rvdsamp{m}=(\rvd_1^{(m)},\ldots,\rvd_d^{(m)})\), and \(w=(w_1,\ldots,w_\ngs)\) be a set of quadrature samples and weights satisfying
\begin{align}\label{eq:gs-quad-rule}
  \int_\dom f(z) \pbwt(z) \dx{z} \approx \sum_{j=1}^\ngs w_j f(\rvdsamp{j}),
\end{align}
and let \(\Psi\) denote the $\ngs \times N$ Vandermonde-type matrix
\begin{align*}
  (\Psi(\nodes_\ngs, N))_{j,n} = \Psi_{j,n} &= \psi_n(\rvdsamp{j}), & j \in [\ngs], &\, n \in [N].
\end{align*}
\rev{Note that designing a quadrature rule with small $J$ so that \eqref{eq:gs-quad-rule} holds is a very difficult task in general. Depending on the approach used, we expect the requisite $J$ to scale with dimension $d$ and with the desired accuracy of \eqref{eq:gs-quad-rule}.
}
We assume that the quadrature rule $(\nodes_\ngs, w)$ forms a proper discrete $\ell^2$ norm on the span of $\psi_n$, or equivalently
\begin{align}\label{eq:orthogonalization-assumption}
  \forall\;\; p \in \mathrm{span}\{ \psi_n\}_{n=1}^N\backslash\{0\}, \hskip 15pt \sum_{j=1}^\ngs w_j p^2(\rvdsamp{j}) > 0.
\end{align}
Under this assumption, the columns of $\Psi$ weighted by $w$ are linearly independent, and so there is a unique QR factorization
\begin{align}
\label{eq:apc-moment-matrix}
\sqrt{W}\Psi=QR,& & W=\text{diag}(w),
\end{align}
\rev{which is effectively performing the operations in \eqref{eq:gram-schmidt-alg}.}
In the remainder of the paper, we will refer to $\sqrt{W}\Psi=QR$ as the moment matrix.
Under these assumptions
\begin{align}
\label{eq:apc-basis}
\phi_n(z)=\sum_{k=1}^n \psi_k(z)(R^{-1})_{kn}, \quad n\in[N]
\end{align}
defines a polynomial basis orthonormal \rev{with respect to the discrete inner product defined by the quadrature rule in \eqref{eq:gs-quad-rule}. For large $J$, this approximates the  arbitrary measure $\pbwt$.} The polynomials $\phi_n$ have the following properties, which follow from \eqref{eq:orthogonalization-assumption}, \eqref{eq:apc-moment-matrix}, and \eqref{eq:apc-basis}:
\begin{itemize}
  \item $\mathrm{span} \{\phi_n\}_{n=1}^P = \mathrm{span} \{\psi_n\}_{n=1}^P$ for any $P \leq N$.
  \item If there is a multi-index set $\Lambda$ such that $\{\psi_n\}_{n=1}^N$ is a basis for $\pi_\Lambda$, then $\{\phi_n\}_{n=1}^N$ is also a basis for $\pi_\Lambda$.
  \item $\sum_{j=1}^\ngs w_j \phi_n(\rvdsamp{j}) \phi_m(\rvdsamp{j}) = \delta_{m,n}$.
\end{itemize}
\noindent From \eqref{eq:apc-basis} it is clear that the basis $\phi_n$ depends on the ordering of the polynomials \(\{\psi_n\}_{n=1}^N\). In the univariate setting, the technique described above is referred to as the Stieltjes procedure \cite{Gautschi_AN_1996}.

\rev{The numerical conditioning of the aforementioned Gram-Schmidt procedure is dependent on the accuracy of the quadrature rule used in~\eqref{eq:gs-quad-rule}. In low dimensions when $\dom$ is tensorial and when the density \(\pbwt\) is known explicitly, we can use tensor-product Gauss quadrature rules to form $\nodes_\ngs$ and $w$. Specifically, to compute a quadrature rule for a dependent measure, we construct a tensor-product quadrature rule for an independent measure $\nu$ that dominates the multivariate dependent measure $\omega$\footnote{By dominate we mean that the dependent measure is absolutely continuous with respect to the tensor-product measure.}. Given the tensor-product rule with points and weights $\{z^{(q)},v_q\}_{q=1}^Q$, we can compute the weighted $L^2_\omega$ inner product necessary to orthogonalize the tensor product basis by a change of variable such that
  \begin{align}\label{eq:gs-tp-quadrature}
    (\psi_n,\phi_k)_{L^2_\omega}&=\int_\Omega \psi_n(z)\phi_k(z) d\omega(z)=\int_\Omega \psi_n(z)\phi_k(z)\omega(z) d\nu(z)\approx \sum_{q=1}^Q{\psi_n(z^{(q)})\phi_k(z^{(q)})}\omega(z^{(q)})v_q\end{align}
    The use of a tensor-product quadrature rule is of course inefficient in general, however it only requires evaluations of polynomials in this context. 

In higher-dimensions, one can use custom polynomial quadrature rules (e.g. ~\cite{Jakeman_N_CMAME_2018}) in the Gram-Schmidt procedure. However, often in these high dimensional settings or when Bayesian inference is used to condition prior estimates of uncertainty on observational data, the joint density of $\rv$ is not known explicitly, but is instead characterized by samples drawn from the unknown distribution. In these situations we must use Monte Carlo or Psuedo Monte Carlo quadrature, where in~\eqref{eq:gs-quad-rule} we use \(w_i=1/\ngs\).}

\section{Sampling schemes for dependent probability measures}
\label{sec-3}
Numerous methods exist for computing the coefficients of a polynomial chaos expansion. In this paper we focus on interpolation-based collocation methods.
In the collocation setting, we have $N = M$, where $N$ is the dimension of the approximation space and $M$ is the number of samples. Given a set of \(N\) realizations  \(\nodes_N=\{\rvdsamp{1},\ldots,\rvdsamp{N}\}\), with 
corresponding model outputs \(y =(f(\rvdsamp{1}),\ldots,f(\rvdsamp{N}))^T\), we would like to find a solution that satisfies \(\Phi\alpha = y\), 
where \(\alpha= (\alpha_{1},\ldots,\alpha_{N})^T\) denotes the vector of PCE coefficients and \(\Phi\in\mathbb{R}^{N\times N}\) denotes the Vandermonde matrix with entries
\(\Phi_{mn} = \phi_n(\rvdsamp{m}),\quad m\in[N],\; n\in[N]\).

Given data \(y\) on \(\nodes_N\), we wish to construct a unique polynomial interpolant \(f_N(\brv)\) from \(\pi_{\Lambda}\) that interpolates \(y\). We assume unisolvence of the interpolation problem for $\nodes_N$ on $\pi_\Lambda$; we will later prescribe a method that numerically guarantees this condition. Assuming unisolvence, there is a unique solution to the linear system:
\begin{align}
  \label{eq:interpolant}
  \Phi \alpha = y, \hskip 10pt \Longleftrightarrow \hskip 10pt V \Phi \alpha = V y, & & V_{ii}=v(\rvdsamp{ii}),i\in[N]
\end{align}
where the latter expression is true when $v(\cdot)$ is a non-vanishing weight function on $\nodes_N$. 
The linear algebraic formulation above can be used to directly motivate a sampling scheme: build $\nodes_M$ sequentially so that the determinant of $V \Phi$ is maximized. \rev{Such a procedure is a particular kind of D-optimal design \cite{fedorov_theory_1972}, but differs from standard D-optimal design constructions in that (i) the procedure is sequential, i.e., greedy maximization is performed so that a \textit{sequence}, instead of a non-nested grid, is constructed, and (ii) the choice of weight $v$ determines special properties of the sequence, as discussed in the following section.} The reason we formulate the interpolation problem \eqref{eq:interpolant} with the additional weight $v$ is both to improve numerical stability, and to form connections with theoretical asymptotic results. For example, the results in \cite{Narayan_J_SISC_2014} show that choosing $v = \sqrt{\pbwt}$ for the construction of weighted Leja sequences produces nodal sets that are asymptotically optimal. 

\subsection{Weighted Leja sequences via LU factorization}
\label{sec-3-2}
A Leja sequence (LS) is essentially a doubly-greedy computation of a determinant maximization procedure. Given an existing set of nodes $\nodes_M$, a Leja sequence update chooses a new node $\rvdsamp{M+1}$ by maximizing the determinant of a new Vandermonde-like matrix with an additional row and column: the additional column is formed by adding a single predetermined new basis element, $\phi_{M+1}$, and the additional row is defined by the newly added point. Hence a LS is both greedy in the chosen interpolation points, and also assumes some \textit{a priori} ordering of the basis elements. The introduction of a row-based weighting makes this a weighted Leja sequence; as the previous section suggests, our weight will be the function $v$; we leave this function undefined for now, but make an explicit choice in Section \ref{sec-3-4}.

In one dimension, a weighted LS can be understood without linear algebra: Let \(\nodes_N\) be a set of nodes on \(\dom\) with cardinality \(N \geq 1\). We will add a new point \(z^{(N+1)}\) to \(\nodes\) determined by the following:
\begin{align}
\label{eq:weighted-leja}
\rvdsamp{N+1} = 
\argmax_{\rvd \in \dom} v(\rvd)\prod_{n=1}^N |\rvd - \rvdsamp{n}|
\end{align}
We omit notation indicating the dependence of $z^{N+1}$ on \(\nodes_N\). 
By iterating \eqref{eq:weighted-leja}, one progressively builds up the Leja sequence \(\nodes\) by recomputing and maximizing the objective function for increasing \(N\). The literature contains many instances of using the above objective function to add nodes \cite{leja_sur_1957,reichel_newton_1990,baglama_fast_1998}.

In multiple dimensions, formulating a generalization of the univariate procedure is challenging. The following linear algebra formulation \cite{Sommiriva_V_CMA_2009,bos_computing_2010} greedily maximizes the weighted Vandermonde-like determinant
\begin{align*}
  \rvdsamp{N+1} = \argmax_{\rvd \in \dom} |\det v(\rvd) \Phi(\nodes, \rvdsamp{N+1})|.
\end{align*}
The above procedure is an optimization with no known explicit solution, so constructing a Leja sequence is challenging. In \cite{Narayan_J_2018}, gradient based optimization was used to construct weighted Leja sequences. However a simpler procedure based upon LU factorization can also be used \cite{bos_computing_2010}. The simpler approach comes at a cost of slight degradation in the achieved determinant of the LS. We adopt the LU-based approach here due to its ease of implementation. 

The algorithm for generating weighted Leja sequences using LU factorization is outlined in Algorithm \ref{alg:lu-leja}. The algorithm consists of 5 steps. First a polynomial basis must be specified. The number of polynomial basis elements must be greater than or equal to the number of desired samples in the Leja sequence, i.e. $N \geq M$. The input basis must also be ordered, and the Leja sequence is dependent on this ordering. Unless otherwise specified, in this paper we only consider total-degree polynomial spaces, that is we have 
\begin{align*}
  \mathrm{span}\{\phi_n\}_{n=1}^N &= \pi_\Lambda, & \Lambda = \Lambda_{k,1}^d,
\end{align*}
for some polynomial degree $k$. We use lexigraphical ordering on $\Lambda$ to define the basis. The second step consists of generating a set of \(S\) candidate samples \(\nodes_S\); ideally, $S \gg M$. Our candidate samples will be generated as independent and identically-distributed realizations of a random variable. The precise choice of the random draw will be discussed in the next section. For now we only require that the measure of the draw have support identical with the measure of $Z$. Once candidates have been generated we then form the $S \times N$ Vandermonde-like matrix \(\Phi\), precondition this matrix with $V$, and compute a truncated LU factorization. (Computing the full LU factorization is expensive and unnecessary.) We terminate the LU factorization algorithm after computing the first \(M\) pivots. These ordered pivots correspond to indices in the candidate samples that will make up the Leja sequence. If we assume that there is \textit{any} size-$M$ subset of $\nodes_S$ that is unisolvent for interpolation, then by the pivoting procedure, a Leja sequence is always chosen so that the interpolation problem is unisolvent. 
\begin{algorithm}
  \caption{Approximate LU sequence}
\label{alg:lu-leja}
\begin{algorithmic}[1]
\REQUIRE number of desired samples \(M\), preconditioning function \(v(\rvd)\), basis \(\{\phi\}_{n=1}^N\)
\STATE Choose the index set \(\Lambda\) such that \(N\ge M\)
\STATE Specifying an ordering of the basis \(\phi\)
\STATE Generate set of \(S\gg M\) candidate samples \(\nodes_S\)
\STATE Build \(\Phi\), \(\Phi_{m,n} =\phi_n(\rvdsamp{m})\), \(m\in[S]\), \(n\in[N]\)
\STATE Compute preconditioning matrix \(V\), \(V_{mm}=v(\rvdsamp{m})\)
\STATE Compute first M pivots of LU factorization, \(PLU=LU(V \Phi\),M) 
\end{algorithmic}
\end{algorithm}
Once a Leja sequence \(\nodes_M\) has been generated one can easily generate a polynomial interpolant with two simple steps. The first step evaluates the function at the samples in the sequence, i.e. \(y=f(\nodes)\). The coefficients of the PCE interpolant can then be computed via
$$\alpha=(LU)^{-1}P^{-1} V y,$$
where the matrices $P$, $L$, and $U$ are identified in Algorithm \ref{alg:lu-leja}.

We end this section by noting that (approximate) Fekete points are an alternative determinant-maximizing choice for interpolation points \cite{Sommiriva_V_CMA_2009,bos_computing_2010,Bos_CLSV_MOC_2011}. We opt to use Leja sequences here because they are indeed a \textit{sequence}, whereas a Fekete point construction is not hierarchical.

\subsection{The induced measure}
\label{sec-3-3}
Generating Leja sequences with large determinants using LU factorization requires generating a large number \(S\) of candidate samples \(\nodes_S\). The only theoretical requirement on the distribution of these samples is that they are sampled over the domain of the random variables \(\dom\). When sampling directly from $\pbwt$ is not feasible, one could instead sample uniformly over $\dom$. However, when a significant portion of the probability is concentrated within a small region of $\dom$, then \(S\) needs to be prohibitively large. Computationally, this manifests as an ill-conditioned matrix $\Phi$. Generation of Leja sequences does not require the evaluation of the expensive simulation model $f$, so that in principle candidate sets can be quite large. However, we wish to avoid LU factorizations of matrices with millions of rows. To reduce the number of candidate samples needed, we instead propose to sample from specific measures. 

Let the samples $\rvdsamp{s}$ be generated as iid realizations from a density $\nu$. In the context of maximizing determinants, an attractive choice of the measure \(\nu\) for generating the samples is the so-called induced measure
\begin{align} 
\label{eq:induced-measure}
\nu(\rvd)=\pbwt(\rvd)k(\rvd) & &k(\rvd)=\sum_{\lambda\in\Lambda} \phi^2_\lambda(\rvd),
\end{align}
where \(k(\rvd)\) is known as the Christoffel function. This measure is a property only of $\pbwt$ and $\pi_\Lambda$, and not of the individual basis elements $\phi_\lambda$. Therefore, this biased density is well-defined even when $Z$ has dependent components. This measure has been shown to generate well-conditioned matrices $\Phi$ \cite{cohen_optimal_2017}. 

Unfortunately, generating samples from the so called induced measure \(v\) requires explicit knowledge of the probability measure \(\pbwt\) of the random variables. In some settings, only samples from the measure are available. To mitigate this issue in this paper, we recommend sampling from a large $\Lambda$-asymptotic distribution called the equilibrium measure. When $\Lambda$ is a total-degree space, the induced measure converges to the equilibrium measure as \(N\rightarrow\infty\) \cite{Narayan_JZ_MC_2017}. In contrast to the induced measure, a simple closed form of the equilibrium measure is known for some distributions. For random variables with bounded densities $\pbwt$ on a hypercube, the equilibrium measure is the Chebyshev distribution. Unfortunately when the probability density $\pbwt$ has areas of high concentration in a hypercube, even generating samples from the equilibrium measure can fail to generate a good candidate set. This is because the equilibrium measure is asymptotically optimal, but may not be actually optimal for finite $\Lambda$. However, we found that enriching the equilibrium set with samples from the probability distribution $\pbwt$ produced excellent candidate sets, and we adopted this approach for all numerical examples in our paper.

\begin{figure}[htb]
\centering
\includegraphics[width=\textwidth]{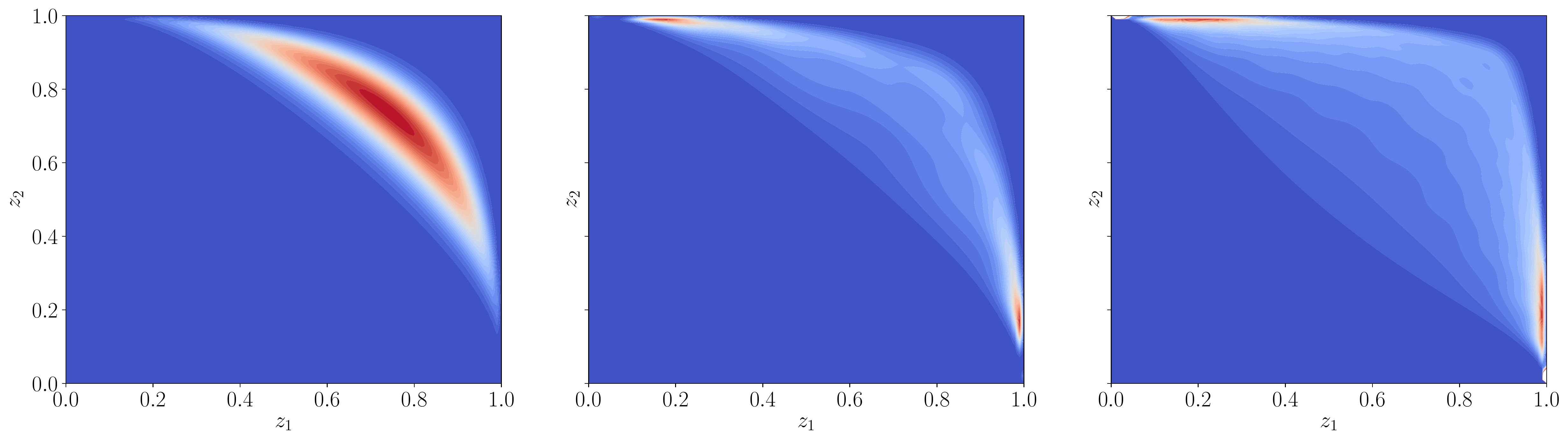}
\caption{(Left) The joint PDF \(\pbwt\) of two variables, given by \eqref{eq:nataf-density}, with \(R^V_{11}=R^V_{22}=1\) and \(R^V_{ij}=-0.9, i\neq j\) and Beta$(2,5)$ marginals. The induced measure density $\nu$ for degree-3 (middle) and degree-20 (right) total degree spaces are also shown.}
\label{fig:cor-beta-2d-density-comparison}
\end{figure}

Let the joint density of $Z$ be given by \eqref{eq:nataf-density}, where the marginal distributions \(\pbwt_i\) are each univariate Beta random variables with parameter $(\alpha,\beta) = (2,5)$ and the correlation matrix  \rev{\(R^V\) with \(R^V_{11}=R^V_{22}=1\) and \(R^V_{ij}=-0.9, i\neq j\)}. The resulting PDF \(\pbwt\) and the induced measure density $\omega_{\Lambda}$ for $\Lambda$ corresponding to degree-3 and degree-20 total degree spaces are depicted in Figure \ref{fig:cor-beta-2d-density-comparison}.

The induced measure assigns significant non-zero probability in regions where the probability of the random variables \(\pbwt\) is small, and the difference between the joint and induced densities appears to increases with polynomial degree. This is consistent with theory that states that for variables on bounded convex domains, the induced distribution will converge (with degree) to a Chebyshev-like distribution \cite{burns_monge-ampere_2010}, which concentrates samples on the boundary of the variable domain \(\dom\). 


\subsection{Weight function}
\label{sec-3-4}
The properties of a Leja sequence are greatly influenced by the choice of the weight function \(v(\brv)\). In \cite{Narayan_J_2018} Leja sequences were generated by setting \(v=\sqrt{\pbwt}\). However in this paper we use the root inverse of the Christoffel function:
\begin{align} 
\label{eq:christoffel-function}
v(\rvd) &= \frac{1}{\sqrt{k(\rvd)}} = \frac{1}{\sqrt{\sum_{\lambda\in\Lambda} \phi^2_\lambda(\brv)}}.
\end{align} 
It is shown in \cite{Guo_NYZ_2018} that the Christoffel function is the optimal choice of weight function for generating Fekete nodes in one dimension. 

\section{Numerical Results}
\label{sec-5}
In this section we demonstrate the efficacy of our proposed approach on a number of numerical examples. 

To measure the performance of an approximation, we will use the \(\pbwt\)-weighted \(\ell^2\) error on a set of test nodes.  We generate a set of \(S=10,000\) random samples \(\{\rvdsamp{j}\}_{j=1}^S\subset\dom\) drawn from the density \(\pbwt\). The error is computed as 
\[
  \label{eq:error}
\lVert f-f_N\rVert_{\ell^2_\pbwt} = \left(\frac{1}{S}\sum_{j=1}^S \lvert f(z^{(j)})-f_N(z^{(j)})\rvert^2\right)^{1/2},
\]
where $f$ is the exact function and $f_N$ is the interpolative approximation. 

When considering dependent random variables with probability measures concentrated in a small region of the variable domain, we are careful not to generate misleading results. Functions that vary significantly in regions of high-probability are much harder to approximate in a weighted \(L^2_\pbwt\) norm than functions that vary in regions of low probability. Here we consider the oscillatory Genz function
\begin{align}
\label{eq:genz-function}
f(\rvd)=\cos(2 \pi \rev{e} + \sum_{i=1}^d c_iz_i), & &z\in[0,1]^d.
\end{align}
This function has strong variation throughout the domain. However, in a further attempt to avoid constructing a function that only varies strongly in regions of low-probability, we set the coefficients \(c\) and \(d\) randomly. Specifically we draw \rev{\(e\)} and \(c_i,i\in[d]\)  randomly from the uniform distribution on \([0,1]\) and set 
\begin{align*}
c_i=\frac{40}{d\sum_{i=1}^d b_i}b_i, & &i\in[d].
\end{align*}

To facilitate notation, let $B(\cdot; \alpha,\beta)$ denote a tensor-product density function of a Beta random variable with parameters $\alpha$ and $\beta$:
\begin{align*}
  B(z;\alpha,\beta) &= K \prod_{i=1}^d z_i^{\alpha-1} (1-z_i)^{\beta-1}, & z &\in[0,1]^d, \;\; \alpha,\beta > 0,
\end{align*}
where $K$ is a normalizing constant to ensure that $B$ is a probability density. Note that we use the same parameters $\alpha, \beta$ for each dimension.

We will consider three approaches in our numerical results:
\begin{itemize}
  \item \texttt{Nataf} -- This is the approach outlined in \eqref{eq:pce-u-space} which uses the transformation $\mathcal{T}_\text{nataf}^\text{gauss}$ defined in Section \ref{sec-4-2}. The approximation $g_N$ from \eqref{eq:pce-u-space} in $U$-space is constructed using weighted Leja sequences as discussed in Section \ref{sec-3}.
  \item \texttt{DOM $(\alpha,\beta)$} -- This is a domination strategy, where the dominating density is $B(\cdot; \alpha,\beta)$. Again, we build approximations with respect to the dominating methods using weighted Leja sequences.
  \item \texttt{GS $(\alpha,\beta)$} -- This is the proposed strategy in this paper. We first construct a basis that is (approximately) $L^2_{\pbwt}$ orthonormal via the technique in Section \ref{ssec:pce-orthogonalization}. We then use weighted Leja sequences as described in Section \ref{sec-3} to construct a PCE approximation in this basis.
\end{itemize}
Unless otherwise stated we use $10,000$ candidate samples to build the Leja sequence, where half of these samples are from a tensor-product Chebyshev density, and the remaining half are drawn from the density $\pbwt$.

We note that all three approches use weighted Leja sequences to construct the approximation. Therefore, our examples are a direct comparison for our three strategies to handle dependent variables: mapping methods, domination methods, and orthogonalization methods.

\subsection{Leja sequence for domination measures}
\label{sec-5-1}
In this section we investigate the impact of constructing a polynomial approximation using the dominating measure strategy of Section \ref{ssec:pce-domination}. Thus given the measure associated with the random variables, we first identify a dominating measure, construct a PCE approximation using an interpolation sequence and basis elements from that dominating measure, and then compute the error with respect to the random variable measure. \rev{Although this paper is an exposition on approximation strategies for dependent random variables, here we choose the measure of the random variables to have independent components to facilitate comparison with known optimal strategies for tensor-product approximation, which are not applicable when variable dependencies exist.} 

\begin{figure}[htb]
\centering
\includegraphics[width=0.46\textwidth]{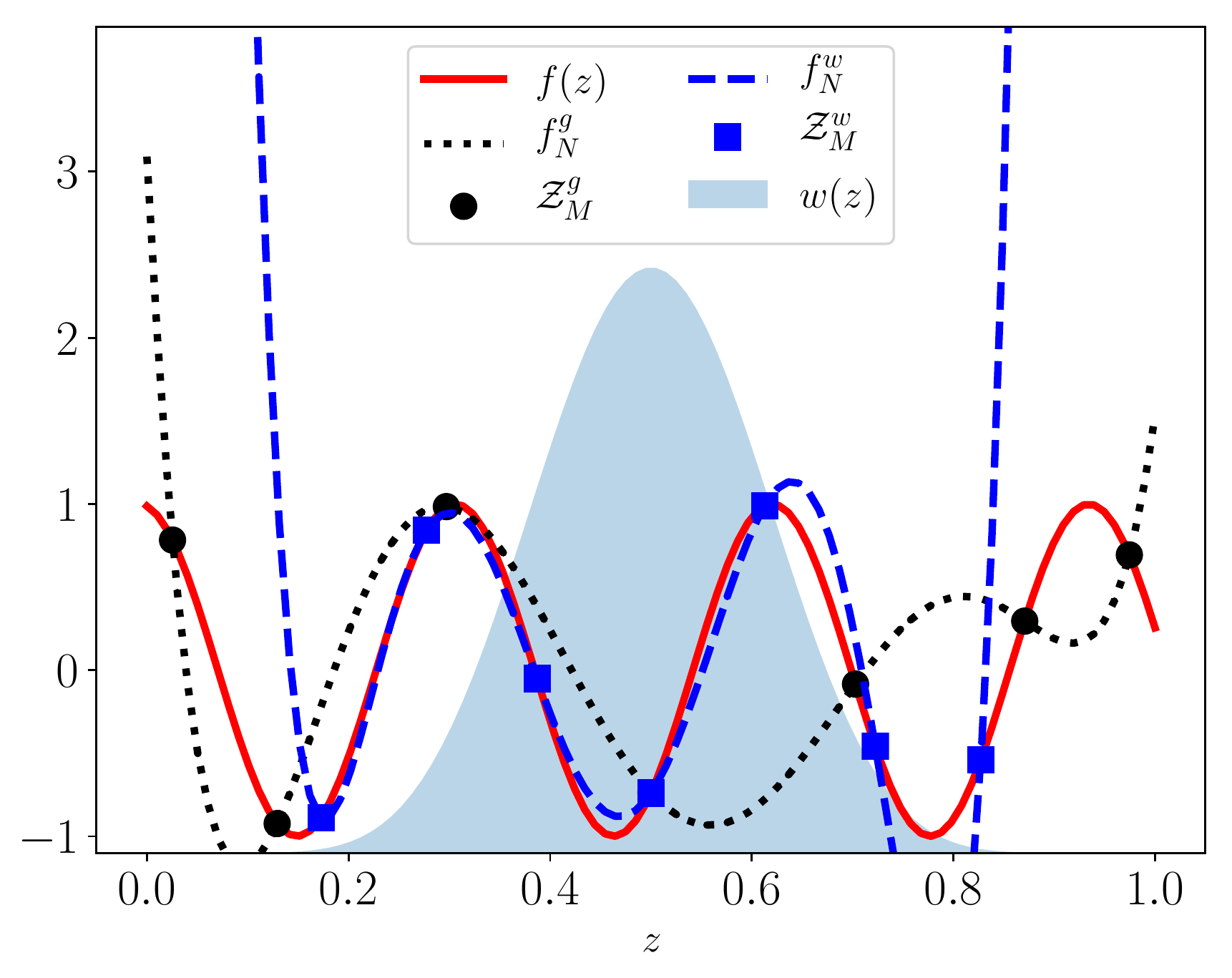}
\includegraphics[width=0.49\textwidth]{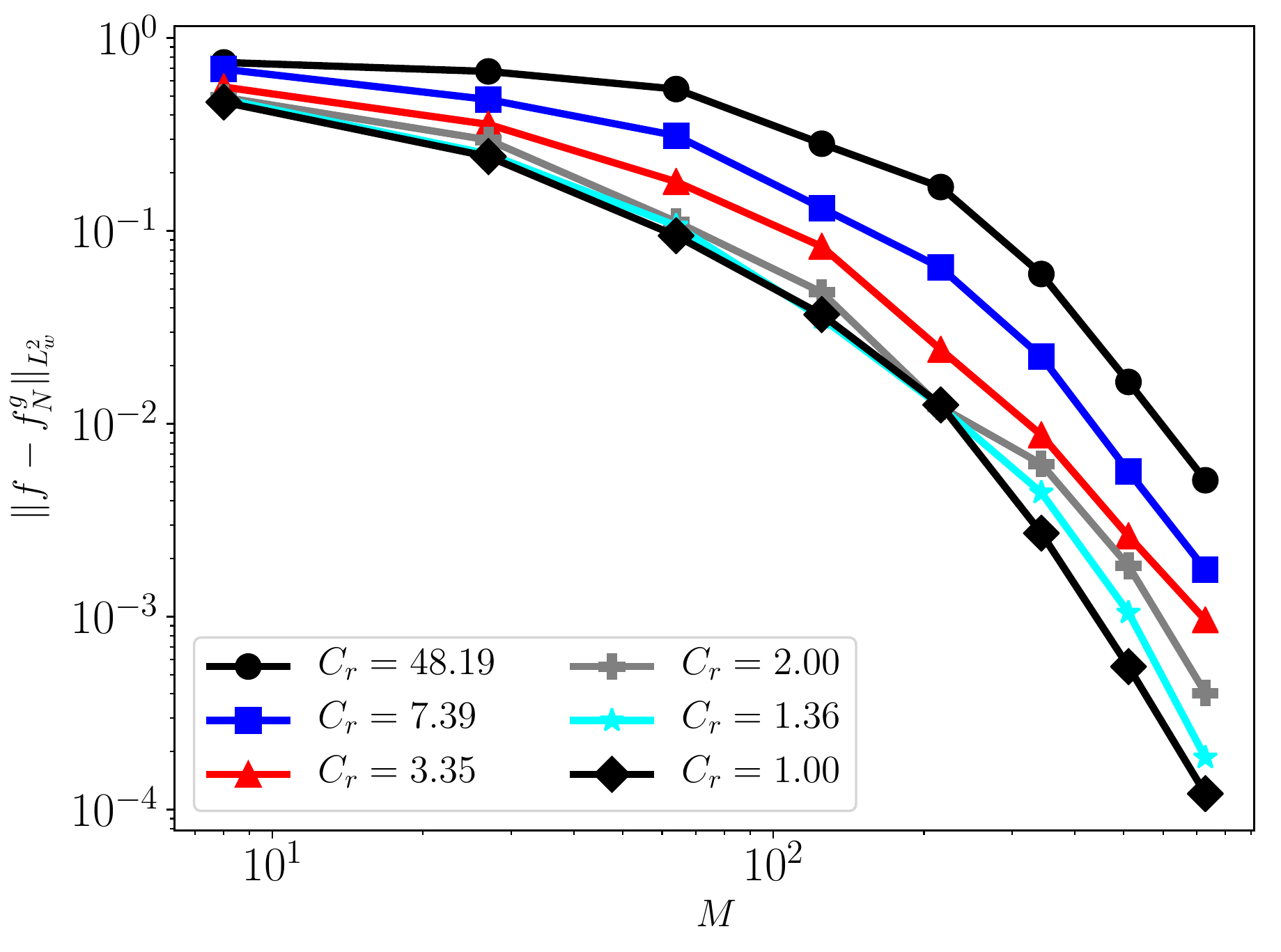}
\caption{(Left) PCE interpolants of a univariate oscillatory Genz function of a Beta random variable with \(\alpha=10, \beta=10\). (Right) Convergence in the error of PCE interpolants based upon varying univariate orthogonal polynomials. In this example $\alpha=\beta$ and the values associated with decreasing $C_r$ are $(0,2,4,6,8,10)$.}
  \label{fig:beta_incorrect_basis_example_1d}
\end{figure}

In Figure \ref{fig:beta_incorrect_basis_example_1d} (left) we consider a one-dimensional case. We compare interpolants of a univariate oscillatory Genz function parameterized by a single Beta random variable $Z$ with the density $\pbwt = B(\cdot; 10, 10)$. Our dominating measure $g$ will be the uniform measure over the same domain, i.e., $g = B(\cdot;1,1)$. We generate two interpolants. The first one, $f_N^\pbwt$, uses Gauss quadrature nodes $\nodes_N^\pbwt$ from $\pbwt$. The second, $f_N^g$, uses Gauss quadrature nodes $\nodes_N^g$ from $g$. However, the error for both is measured with the same formula and samples in \eqref{eq:error} that are generated from the Jacobi measure. The Gauss-Jacobi quadrature samples are optimal in the univariate setting and will produce a unitary condition number. The sampling scheme of the Jacobi basis places samples in a way that balances stability with sampling in high-probability regions. The procedure with $g$ does this as well, but for the uniform density. Consequently the samples for the uniform measure appear more frequently in regions of low $\pbwt$-probability. This manifests itself in a larger error in regions of high $\pbwt$ probability than obtained using the optimal $\pbwt$ construction. It is also apparent from Figure \ref{fig:beta_incorrect_basis_example_1d} (left) that the approximation $f_N^\pbwt$ peforms poorly in regions of very low $\pbwt$ probability. However, as measured in the $L^2_{\pbwt}$ norm, this is allowable.

In Figure \ref{fig:beta_incorrect_basis_example_1d} (right) we also plot the convergence in the \rev{median error (over $10$ samples of $b_i$ and $e$)} of PCE approximations of the analytic function \eqref{eq:genz-function} in three dimensions, $d=3$. We set the joint density as \(\pbwt = B(\cdot;10,10)\), and construct a PCE via the measure domination technique in Section \ref{ssec:pce-domination}, with the dominating density $g = B(\cdot;\beta,\beta)$ set to be another tensor-product Beta measure over the same support.  We use various values of $\beta < 10$. As $\beta$ approaches $10$ the dominating measure approximation will become more efficient. To limit the effect of the sampling scheme we use tensor-products of the univariate Gauss quadrature rule for $g$, and set the approximation space as associated with a tensor-product index set approximation,
$$
\Lambda = \{\phi_{\lambda} : \norm{\lambda}{\infty} \le p\}=\Lambda_{p,\infty}.
$$

For each choice of dominating measure we compute the constant \(C_r\) from Lemma \ref{lem:chen-2013-strong} as a measure of the distance between the measure of orthogonality of each polynomial basis and the probability measure \(\pbwt\). The plot clearly shows that constructing an approximation from a measure that is not orthogonal to \(\pbwt\) results in a degradation of accuracy, and the penalty grows as \(C_r\) grows. Using a dominating measure effects the constant of convergence but not the rate of convergence which is consistent with Lemma~\ref{lem:chen-2013-strong}.

\subsection{Leja sequences using probabilistic transformations}
\label{sec-8-1}
In this section we explore the performance of interpolants constructed using the transformations $\mathcal{T}_\text{nataf}^\text{gauss}$ and $\mathcal{T}_\text{nataf}^\text{unif}$, defined in Section \ref{sec-4-2}, along with the mapping procedure \eqref{eq:pce-u-space} and weighted Leja sequences.

\begin{figure}[htb]
\centering
\includegraphics[width=\textwidth]{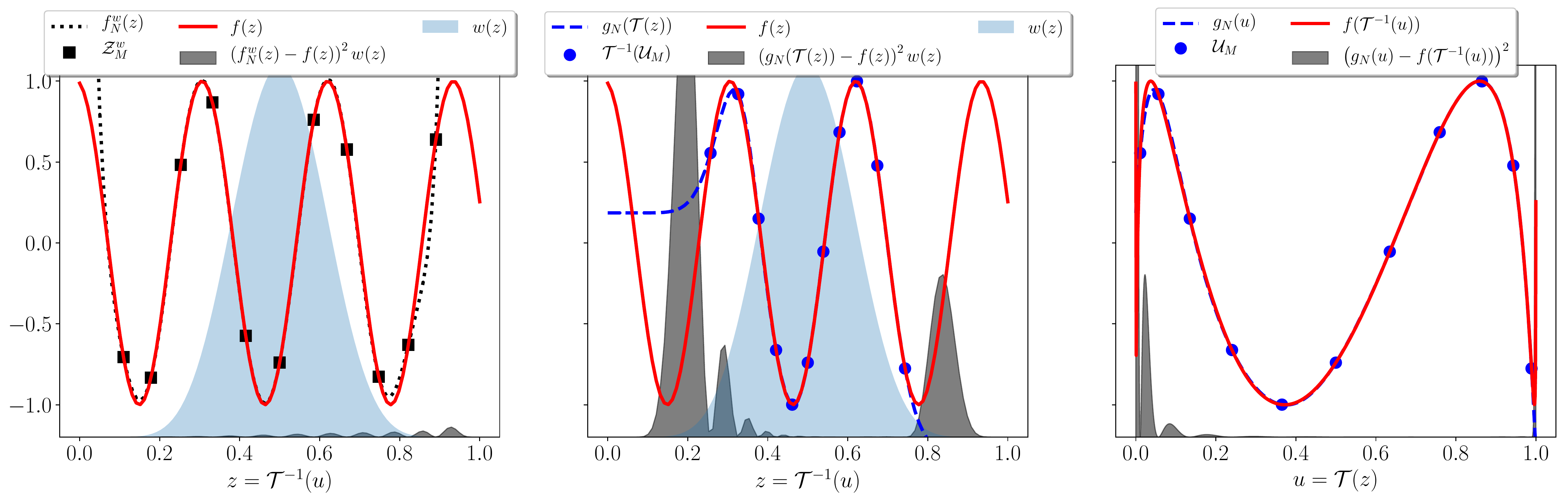}
\caption{Gauss quadrature interpolant of a cosine function (left) versus mapped methods using Gauss quadrature interpolants (center, right). The pointwise error (grey shaded region) of the mapped methods is significantly larger than the direct Gauss quadrature method. This suggests that, for the same number of samples, mapping methods produce suboptimal approximations. }
  \label{fig:rosenblatt_1d_beta_example}
\end{figure}

Probabilistic transformations, such as the Nataf transformation, are highly non-linear and often introduces steep gradients into a function via the map composition. To highlight the effect of the increased non-linearity introduced by the Nataf transform, consider the oscillatory Genz function as a function of only one variable shown in Figure \ref{fig:rosenblatt_1d_beta_example}. In this case the Rosenblatt and $\mathcal{T}_\text{nataf}^\text{unif}$ transformations are equivalent and consist of simply applying inverse transform sampling via the distribution function to define the \(Z\) to \(U\) transformation and back. Here we approximate the function using Gauss quadrature nodes of polynomials orthonormal with respect to the variables \(Z\) and \(U\). These are optimal in one-dimension and will outperform even Fekete and Leja sequences. Thus, we have two approximations: $f_N^{\pbwt}(z)$ which is a Gauss quadrature interpolant built from the Gauss nodes of $\pbwt$; and $g_N(\mathcal{T}(z))$, where $g_N$ is a Gauss quadrature interpolant built from Gauss nodes of a uniform random variable, and subsequently mapped to $z$-space via composition with $\mathcal{T}=\mathcal{T}_\text{nataf}^\text{unif}$.

We see the strong non-linearity introduced by the transformation degrades the accuracy of the interpolant. We can also easily see that the approximation error is largest in the regions in which the non-linearity introduced is strongest. These regions occur in regions of significant probability and so the resulting \(L^2_\pbwt\) error is larger when using the Nataf transformation than without.

\begin{figure}[htb]
\centering
\includegraphics[width=\textwidth]{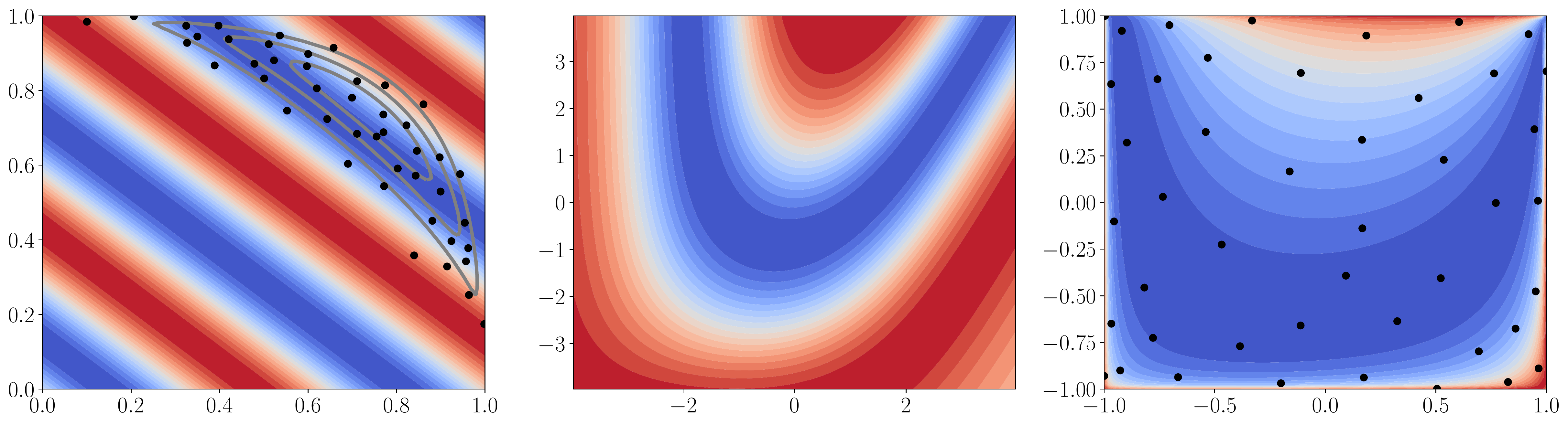} 
\caption{\rev{Comparison of the algebraic function \(f(\brv)\) in \eqref{eq:genz-function} (left) and the same function under the transformations $\mathcal{T}_\text{nataf}^\text{gauss}$ (middle) and $\mathcal{T}_\text{nataf}^\text{unif}$ (right). The contours in the left plot represent the probability density given by \eqref{eq:nataf-density}, with \(R^V_{11}=R^V_{22}=1\) and \(R^V_{ij}=-0.9, i\neq j\) and Beta$(2,5)$ marginals. The dots are the Leja samples in the space of the original dependent variables \(\brv\) (left) and the independent variables \(U\) (right).}}
  \label{fig:nataf-function-comparison}
\end{figure}

Figure \ref{fig:nataf-function-comparison} depicts the two-dimensional version of the algebraic function defined in \eqref{eq:genz-function}. \rev{The same figure also plots the function under the transformations $\mathcal{T}_\text{nataf}^\text{gauss}$ and $\mathcal{T}_\text{nataf}^\text{unif}$. Clearly both transformations} increase the non-linearity of the function, which means that we must also increase the degree of the PCE and thus the number of samples needed to approximate the function accurately.
In Figure \ref{fig:nataf-function-comparison} we also superimpose a degree 12 Leja sequence constructed in the i.i.d uniform space, along with the same samples mapped back to the original space \(Z\) on the right and left plots respectively. This Leja sequence produces an interpolation matrix with a small condition number, however this benefit is outweighed by the degradation in accuracy of the interpolant when compared to the dominating measure approach (see Figure \ref{fig:cor-beta-2d-interpolation-condition-num-comparison}). \rev{Note that we can also build Leja sequence using univariate Hermite polynomials in the i.i.d. standard normal space coupled with the transformation $\mathcal{T}_\text{nataf}^\text{gauss}$. But again due to the non-linearity of the Nataf transformation (see Figure~\ref{fig:nataf-function-comparison}), the accuracy of the interpolant is poor relative to the approximation constructed using the dominating measure approach. We omit this comparison for brevity.}

Another way to interpret the slow convergence of the interpolant built using non-linear probabilistic transformations is that we are no-longer approximating with polynomials in the original space of \(\brv\), but rather with mapped polynomials, and hence do not have the same guarantees on convergence that exist when approximating with polynomials. The poor performance we observe here is consistent with the univariate results documented in \cite{Xiu_K_SISC_2002} and the multivariate results reported in \cite{Eldred_B_AIAA_2009}.

\subsection{Leja sequences using GSO bases}
\label{sec-5-2}
In this section we compare our three approaches for dependent variables when the approximation is constructed using interpolation via weighted Leja sequences. To consider the effect of dimension we consider two problems with \(d=2\) and \(d=10\). In each case we assume that the joint density of $Z$ is given by \eqref{eq:nataf-density} and that the marginal distributions \(\pbwt_i\) are each univariate Beta random variables with parameters $(\alpha,\beta) = (2,5)$.

\begin{figure}[htb]
\centering
\includegraphics[width=\textwidth]{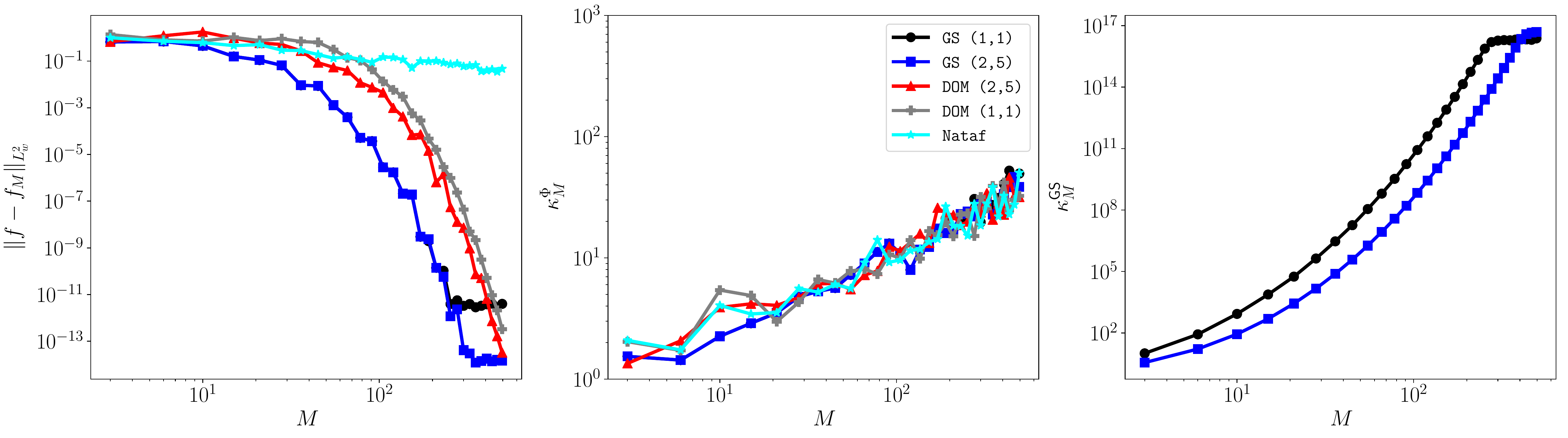}\\
\includegraphics[width=\textwidth]{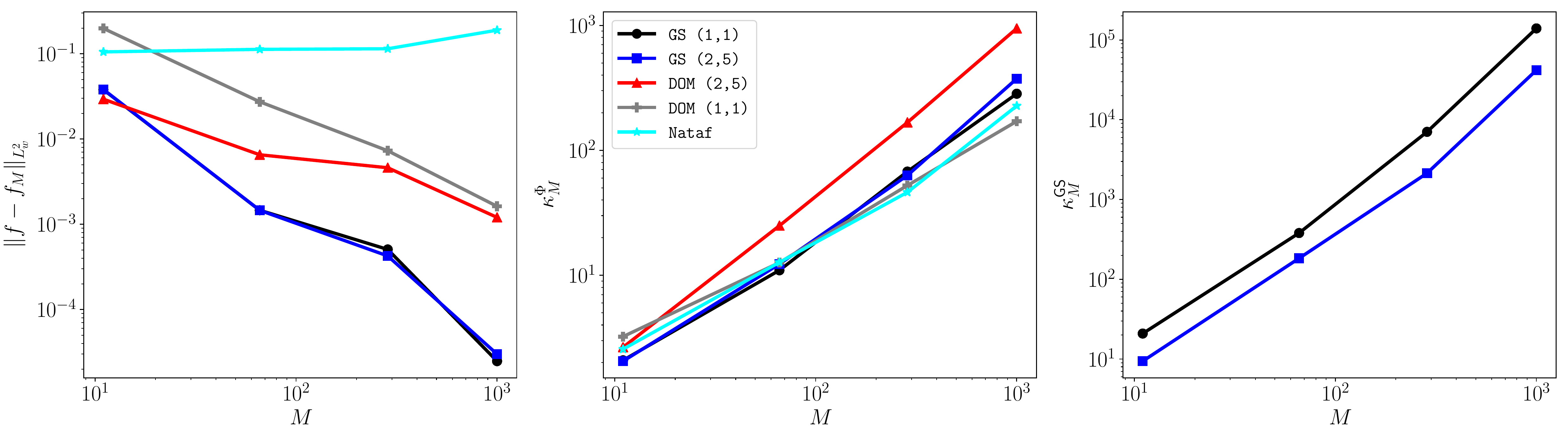}
\caption{(Left) \(L^2_\pbwt\)-errors in the polynomial interpolants of the 2D oscillatory Genz function. (Middle) condition number of the Vandermonde matrix \(\Phi(\nodes_N)\) evaluated at the interpolation points. (Right) for the orthogonalization methods, condition number of the weighted moment matrix $\sqrt{W} \Phi$ in \eqref{eq:apc-moment-matrix}. Each curve represents the median over 10 trials. Top $d=2$, Bottom $d=10$. Recall the notation $(\alpha,\beta)$ used in the legend refers to the parameters of the Beta distribtion that is the orthogonality measure of univariate polynomials used in the GSO procedure.}
\label{fig:cor-beta-2d-interpolation-condition-num-comparison}
\end{figure}

In the top row of Figure \ref{fig:cor-beta-2d-interpolation-condition-num-comparison} we compare the $\texttt{Nataf}$, $\texttt{DOM}$, and $\texttt{GS}$ approaches described earlier for building PCE interpolants using the test function \eqref{eq:genz-function} in two dimensions, \rev{where we set the entries of the correlation matrix \(R^V\) to \(R^V_{11}=R^V_{22}=1\) and \(R^V_{ij}=-0.9, i\neq j\)}. This density is plotted in the left of Figure~\ref{fig:cor-beta-2d-density-comparison} and right of Figure \ref{fig:nataf-function-comparison}.
With the exception of the Nataf interpolant, the rates of convergence in the \rev{median error (over $10$ samples of $b_i$ and $e$)} are similar. However the improved constant of convergence in \eqref{eq:episetmic-theorem-result} obtained by using the GS-orthogonalized basis results in significant improvements in accuracy. The Leja sequences built using the GSO basis produce interpolants that have errors which are orders of magnitude smaller for a fixed sample size. Unlike the other approaches, the convergence rate of the error in the Nataf interpolant deteriorates significantly. The poor performance we observed here is consistent with results observed in the univariate case \cite{Xiu_K_SISC_2002} and the multivariate results observed in \cite{Eldred_B_AIAA_2009}. In short, this is caused by a Jacobian with steep gradients, as explained in Section~\ref{sec-8-1}.

To assess the conditioning of the interpolation problem and the procedure used to construct the GS basis we respectively use the following condition numbers
\begin{align}
\kappa^\Phi_N=\frac{\sigma_{\max}(V\Phi)}{\sigma_{\min}(V\Phi)}, & & \kappa^{\text{GS}}_N=\frac{\sigma_{\max}(\sqrt{W}\Psi)}{\sigma_{\min}(\sqrt{W}\Psi)}
\end{align}
where $\sigma_{\max}(A)$ and $\sigma_{\min}(A)$ are the maximum and minimum singular values of a matrix A.
In all cases the weighted Vandermonde-like matrices $V\Phi$ built using Leja sequences are well conditioned. 
This is despite the fact that the condition numbers $\kappa^{\text{GS}}_N$ of the moment matrices \eqref{eq:apc-moment-matrix} used to compute the GSO basis are large. The condition number of the moment matrices does eventually affect accuracy of the interpolant, but only at extremely high degrees, specifically \(p\ge 23\) for the \texttt{GS $(1,1)$} case. 

The condition number of the moment matrix is dependent on how disparate the orthogonality measure of the tensor product basis is from $\pbwt$. The condition number can be reduced by using a tensor product basis which is orthogonal to a measure which is ``closer'' to the probability measure \(\pbwt\). This is evident by the reduction in condition number in the \texttt{GS $(2,5)$} example. The problem tested here is quite challenging, i.e. the correlation of the variables is high and the ratio of the densities is large, but our proposed algorithm still performs well.

The bottom row of Figure \ref{fig:cor-beta-2d-interpolation-condition-num-comparison} plots the convergence of the various interpolation strategies for the algebraic equation with \(d=10\). Each off diagonal entry of $\Rvy$ is set to 0.9, then for $i=1,3,\ldots,9$, $j=1,\ldots,d$ we set $\Rvy_{ij}=-\Rvy_{ij}$, and $\Rvy_{ji}=-\Rvy_{ji}$. We observe similar trends to those in the 2D example. Specifically the GS approaches perform notably better than the other approaches. Again, the Nataf transformation performs poorly: after 1000 samples the error is still \(\mathcal{O}(1)\). As in the 2D case, the \texttt{GS $(2,5)$} case produces better Gram-Schmidt condition numbers than the \texttt{GS $(1,1)$} case.

\subsection{Using Leja sequences for quadrature}
\label{sec-5-3}
The construction of Leja points is motivated largely by interpolation; however it is straightforward to compute the mean and variance from a polynomial chaos expansion constructed via interpolation. Due to orthonormality, the mean and variance can be computed analytically from the PCE
\begin{align}
\label{eq:pce-moments}
\mu_{f_N}=\mathbb{E}_\brv[f_N]=\alpha_{(0,\ldots,0)} & & \sigma^2_{f_N}=\mathrm{Var}_\brv[f_N]=\sum_{\lambda\in\lambda}\alpha_\lambda^2-\mu^2_{f_N}
\end{align}
\rev{In the case when an approximate grid is used to orthogonalize polynomials, such as in \eqref{eq:gs-quad-rule}, then the expressions above are only true with respect to the discrete measure defined by this grid. If the discrete grid is accurate, then the above can be reasonable approximations to expectations with respect to the continuous dependent density.}
If the PCE interpolates using a set of nodes \(\nodes_N\), then we can also easily compute a quadrature rule. Because we can compute mean and variance analytically using \eqref{eq:pce-moments} such rules are not often needed. However the weights of the quadrature rule can be used to quantify the ill-conditioning of the estimates of integrals using PCE interpolants.

Given function data \(y\), the polynomial chaos coefficients \(\alpha_n\) satisfy
\begin{align*}
  \Phi \alpha &= y, & \Phi_{m,n} = \phi_{n}(\rvdsamp{m}).
\end{align*}
Since \(\alpha = \Phi^{-1} y\), and \(\phi_1 \equiv 1\) because \(\pbwt\) is a probability density function,
\begin{align*}
  \int_{\dom}\sum_{n=1}^{N} \alpha_n \phi_{n}(\rvd) \pbwt(\rvd) \dx{z} = \alpha_1 \int_{\dom} \phi_1(\rvd) \pbwt(\rvd) \dx{z} = \alpha_1,
\end{align*}
then we immediately conclude that the first row of the matrix \(\Phi^{-1}\) gives us quadrature weights \(v_n\) defining the Leja polynomial quadrature rule
\begin{align*}
  \int_{\dom} f_N(\rvd) \pbwt(\rvd) \dx{\rvd} = Q_{N} f = \sum_{n=1}^{N} v_n f(\rvdsamp{n}).
\end{align*}
The subscript \(N\) on $Q_N$ indicates that the quadrature rule has \(N\) nodes. 

Given the weights \(v_n\), the $\ell^1$ condition number of the quadrature operator $Q_N$ is given by 
\begin{align}
\label{eq:quad-rule-condition-number}
\kappa_{Q_N} = \frac{\sum_{n=1}^N \lvert v_n \rvert}{\sum_{n=1}^N v_n}=\sum_{n=1}^N \lvert v_n \rvert
\end{align}
where the last equality holds under the assumption that the \(\phi_n\) are orthonormal with respect to a probability density function \(\pbwt\). Large values of \(\kappa_{Q_N}\) indicate the presence of negative weights, which makes the computation susceptible to catastrophic cancellation.

\begin{figure}[htb]
\centering
\includegraphics[width=\textwidth]{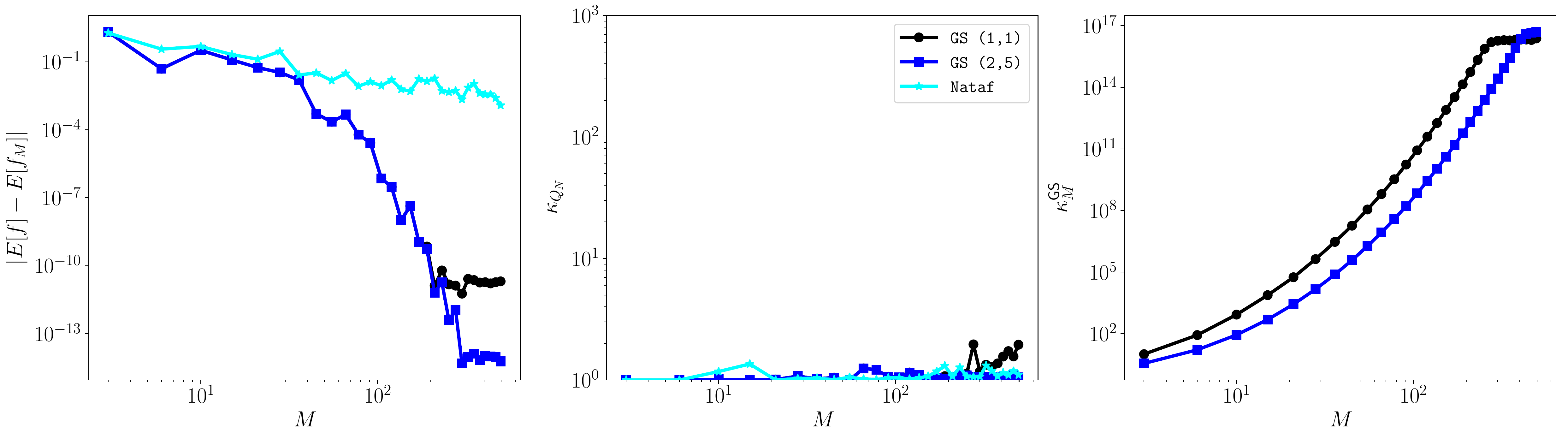}
\includegraphics[width=\textwidth]{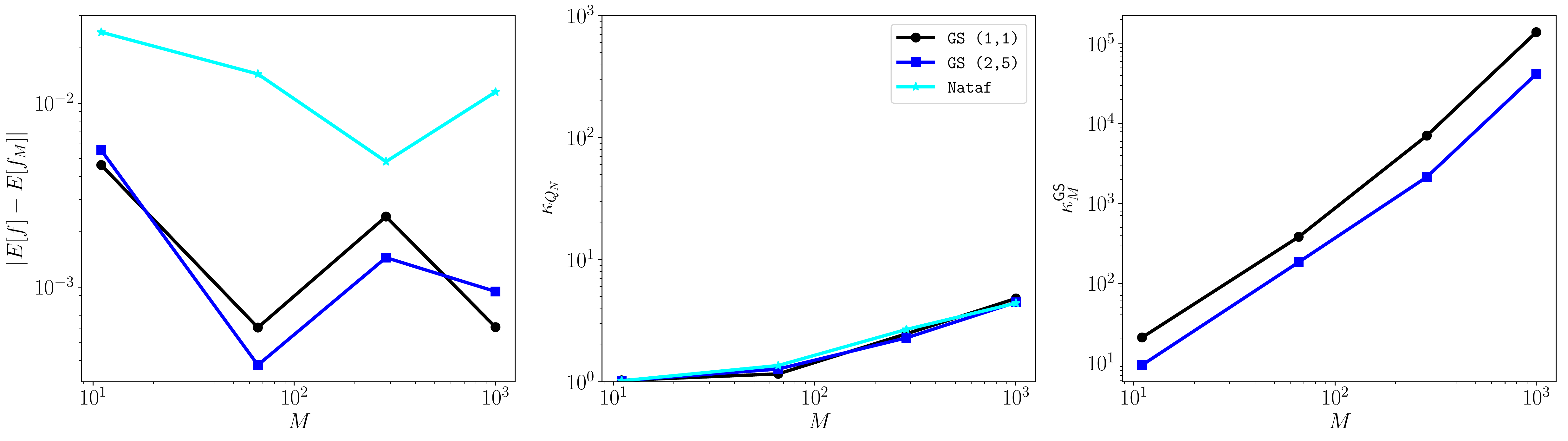}
\caption{(Left) Relative errors in the mean of PCE interpolants of the 2D function \eqref{eq:genz-function}. (Middle) condition number \(\kappa_{Q_N}\) of the quadrature rule. (Right) condition number of the weighted moment matrix $\sqrt{W} \Phi$ from \eqref{eq:apc-moment-matrix}. Each curve represents the median over 10 trials. Top $d=2$, Bottom $d=10$.}
\label{fig:cor-beta-2d-quadrature-condition-num-comparison}
\end{figure}

In Figure \ref{fig:cor-beta-2d-quadrature-condition-num-comparison}, we consider the same examples as in Figure \ref{fig:cor-beta-2d-interpolation-condition-num-comparison}, but now plot the relative error in the mean of the PCE obtained using different Leja sequences and the quadrature rule condition numbers. The \texttt{DOM $(1,1)$} and \texttt{DOM $(2,5)$} approaches cannot directly estimate the mean since the basis elements are not orthogonal to the probability measure, so we do not plot the errors for these approaches (integrals can be estimated numerically by sampling on the interpolant from the correct measure). As observed when computing weighted \(L^2_\pbwt\) errors, the Nataf transformation performs poorly when used for quadrature and the Gram-Schmidt Leja sequences produce the most accurate estimates. The condition number of all quadrature rules are observed to be \(O(1)\). However we must note that again the condition number of the moment matrices used to construct the Gram-Schmidt basis affect accuracy at very high-degrees. In two-dimensions, we compute the true mean of the function to machine precision using Gauss-Legendre tensor product quadrature, that is, we evaluate $\int_\dom f(\rvd) \dx{\pbwt(\rvd)}$ by integrating $f(\rvd) \pbwt(\rvd)$ with respect to the uniform measure on $[0,1]^2$. In ten dimensions, tensor product quadrature is not feasible so we use Monte Carlo (MC) quadrature with $10^6$ samples drawn from the probability measure. The non-monotonic decrease in the error of the Gram-Schmidt PCE rules results from the error in the MC estimate of the mean being larger than the error in the PCE estimates.

\subsection{Approximation without a closed form for the joint density}
\label{sec-5-4}
When building an interpolant of a simulation model, a closed form of the joint density of the random variables is often not available. Instead one may only have samples from the joint density. For example, Bayesian inference \cite{Stuart_AN_2010} is often used to infer densities of random variables conditional on available observational data. The resulting so-called posterior densities are almost never tensor-products of univariate densities. Moreover it is difficult to compute analytical expressions for the posterior density and so Markov Chain Monte Carlo (MCMC) sampling is often used to draw a set of random samples from the posterior density.

\begin{figure}[htb]
\centering
\includegraphics[width=\textwidth]{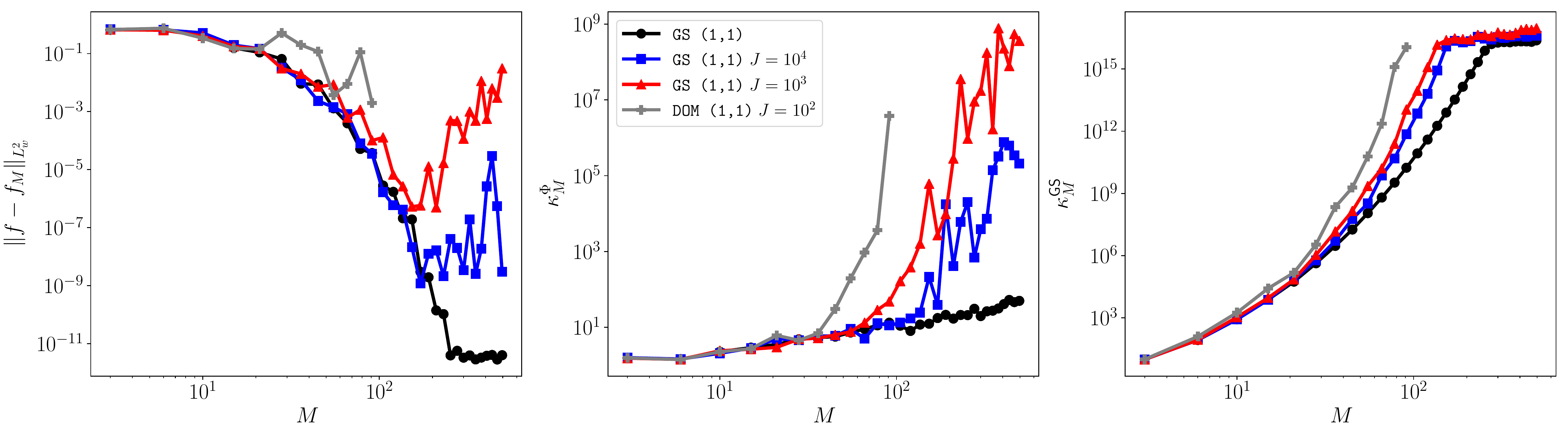}\\
\includegraphics[width=\textwidth]{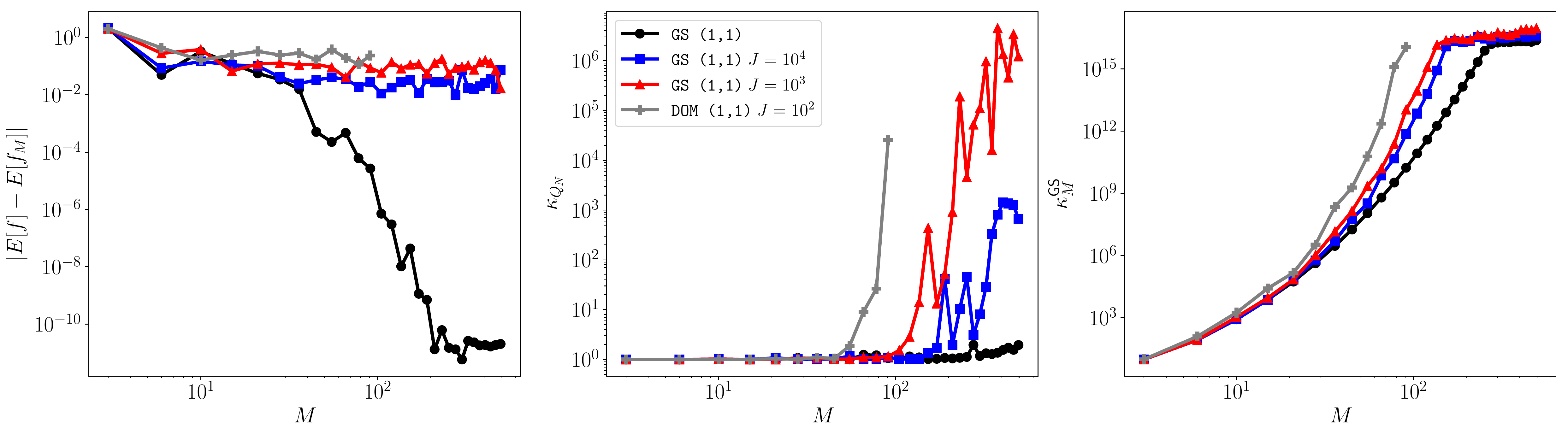}
\caption[]{Metrics for the \texttt{GS $(2,5)$} procedure for approximating the function \eqref{eq:genz-function}. (Top-left) \(L^2_\pbwt\)-errors in the polynomial interpolants. (Top-middle) condition number of the weighted Vandermonde matrix \(V \Phi(\nodes_M)\) evaluated at the interpolation points. 
(Bottom-left) relative errors in the mean of the PCE interpolants. (Bottom-middle) condition number \(\kappa_{Q_N}\) of the quadrature rule.
(Top and bottom right) condition number of the weighted moment matrix in \eqref{eq:apc-moment-matrix}. Each curve represents the median over 10 trials. $\ngs$ refers to the number of Monte Carlo samples used to compute the GS basis. When no $\ngs$ is given \eqref{eq:gs-tp-quadrature} was used.}
\label{fig:mc-sample-size-comparison}
\end{figure}

Here we investigate the impact of computing the moment matrix \eqref{eq:apc-moment-matrix}, used with the Gram-Schmidt orthogonalization, using Monte Carlo (MC) quadrature when no closed form of the joint density is known. In the top row of Figure \ref{fig:mc-sample-size-comparison} we compare the accuracy of the interpolation procedure, when applied to the 2D algebraic function, as the number \(J\) of MC samples is varied. From the figure we can see that decreasing the number of MC samples increases the ill-conditioning of both the  interpolation matrix and the moment matrix. When the condition numbers become sufficiently large the ill-conditioning also affects the accuracy of the interpolant. Clearly the onset of the degradation of the accuracy of the interpolant is dependent on the number of MC samples.

The bottom row of Figure \ref{fig:mc-sample-size-comparison} highlights the effect of the number of MC samples on the estimation of the mean of the 2D algebraic function. The accuracy of the mean obtained from the PCE is limited by the number of samples used to construct the orthonormal basis. This is consistent with the observations made in the univariate setting in \cite{Oladyshkin_N_RESS_2012}. The effect of the sample size appears to have a much greater impact on the accuracy of the quadrature procedure than it does on interpolation accuracy measured in the $L^2_{\pbwt}$ norm. This fact can be leveraged to accurately estimate statistics using the GSO procedure by randomly sampling on the PCE approximation. This sampling can be done with no additional expensive evaluations of the true function $f$. Note the curve depicting approximation with 100 MC samples does not extend to large degrees because a requirement of the Gram-Schmidt orthogonalization is that the moment matrix is over-determined, which is not true for large degrees when only 100 MC samples are used.

\subsection{Approximating a chemical reaction model}
\label{sec-5-5}
In this section we will demonstrate the utility of using our approach for Bayesian inference and dimension reduction using a model of competing species absorbing onto a surface out of a gas phase \cite{Vigil_W_PRE_1996}. Consider the following ordinary differential equation, prescribing evolution of the mass fractions $(u_1, u_2, u_3)$ of three chemical species:
\begin{align}
\label{eq:chemical-species}
\begin{split}
\frac{du_1}{dt} &= as-cu_1 - 4du_1u_2\\
\frac{du_2}{dt} &= 2bs^2 - 4du_1u_2\\
\frac{du_3}{dt} &= es - fu_3\\
s=u_1&+u_2+u_3,
\end{split}
\end{align}
for some constants \(a\), \(b\), \(c\), \(d\), \(e\), and \(f\) and initial conditions \( u_1(0),u_2(0),u_3(0)\).

\subsubsection{Bayesian inference}
\label{sec-5-5-1}
Bayesian inference \cite{Stuart_AN_2010} is often used to infer densities of random variables conditional on available observational data. To make this precise, let \(f_\text{o}({\brv}) : \R^{d}\rightarrow\R^{n_o}\) be an observable quantity, parameterized by the same \(d\) random variables \({\brv}\), which predicts a set of \(n_o\) observable quantities. Bayes rule can be used to define the posterior density for the model parameters \({\brv}\) given observational data \({{y}_o}\):
\begin{align}
\label{eq:posterior}
\pi({\rvd}\mid {{y}_o})=\frac{\pi({{y}_o}\mid {\rvd})\pi({\rvd})}{\int_{\dom}
\pi({{y}_o}\mid {\rvd})\pi({\rvd})d{\rvd}},
\end{align}
where any prior knowledge on the model parameters is captured 
through the prior density \(\pi(\rvd)\). The function \(\pi(y_0 \mid \rvd )\) is the likelihood function and dictates an assumed model-versus-data misfit.

The construction of the posterior is often not the end goal of an analysis. Instead, one is often interested in statistics on the unobservable QoI. Here we will
focus on approximation of the data-informed predictive distribution of the mass fraction of the first species $u_1(t=50)$ at $t=50$ seconds, which we assume we cannot measure. 

\begin{figure}[htb]
\centering
\includegraphics[width=0.48\textwidth]{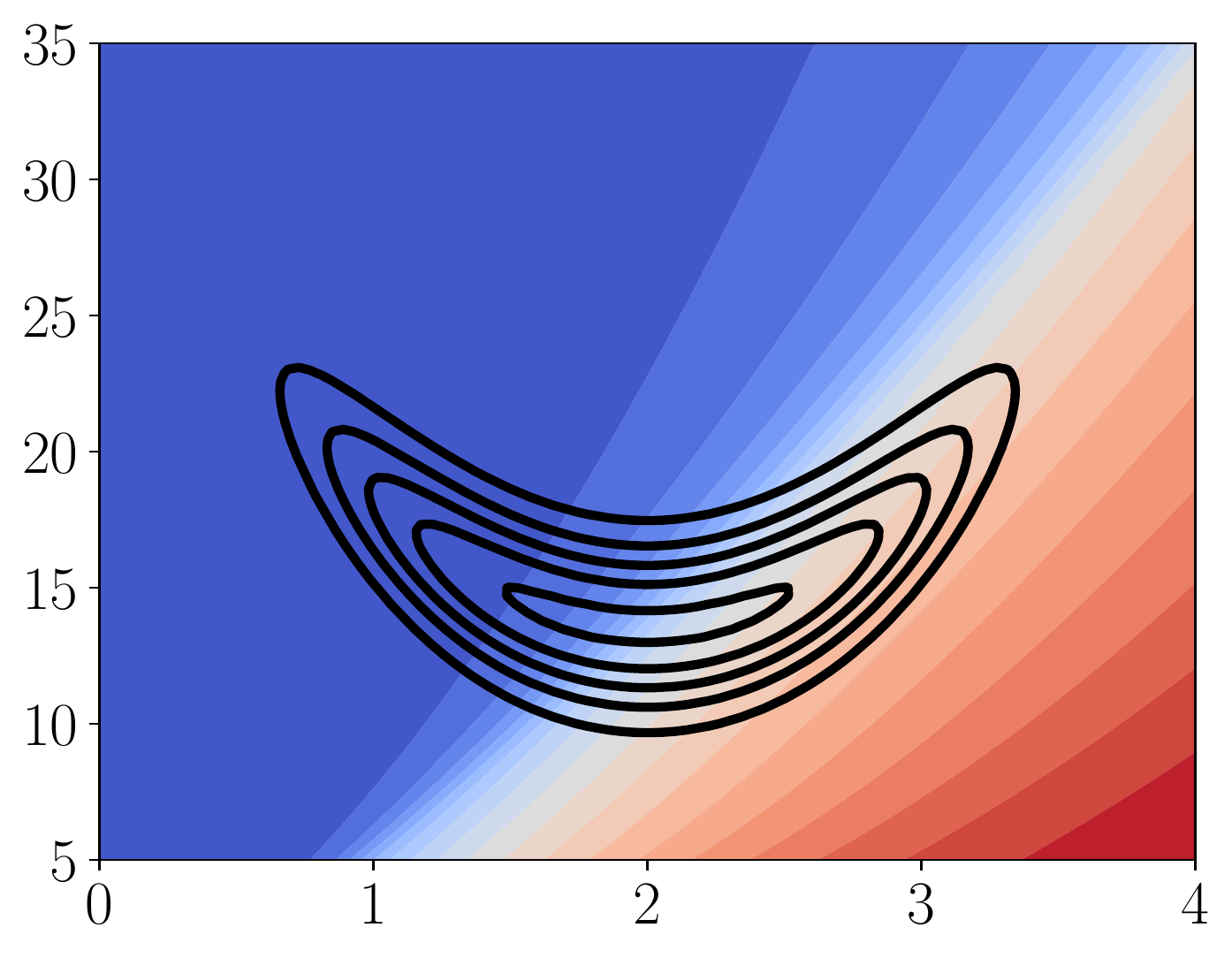}
\includegraphics[width=0.48\textwidth]{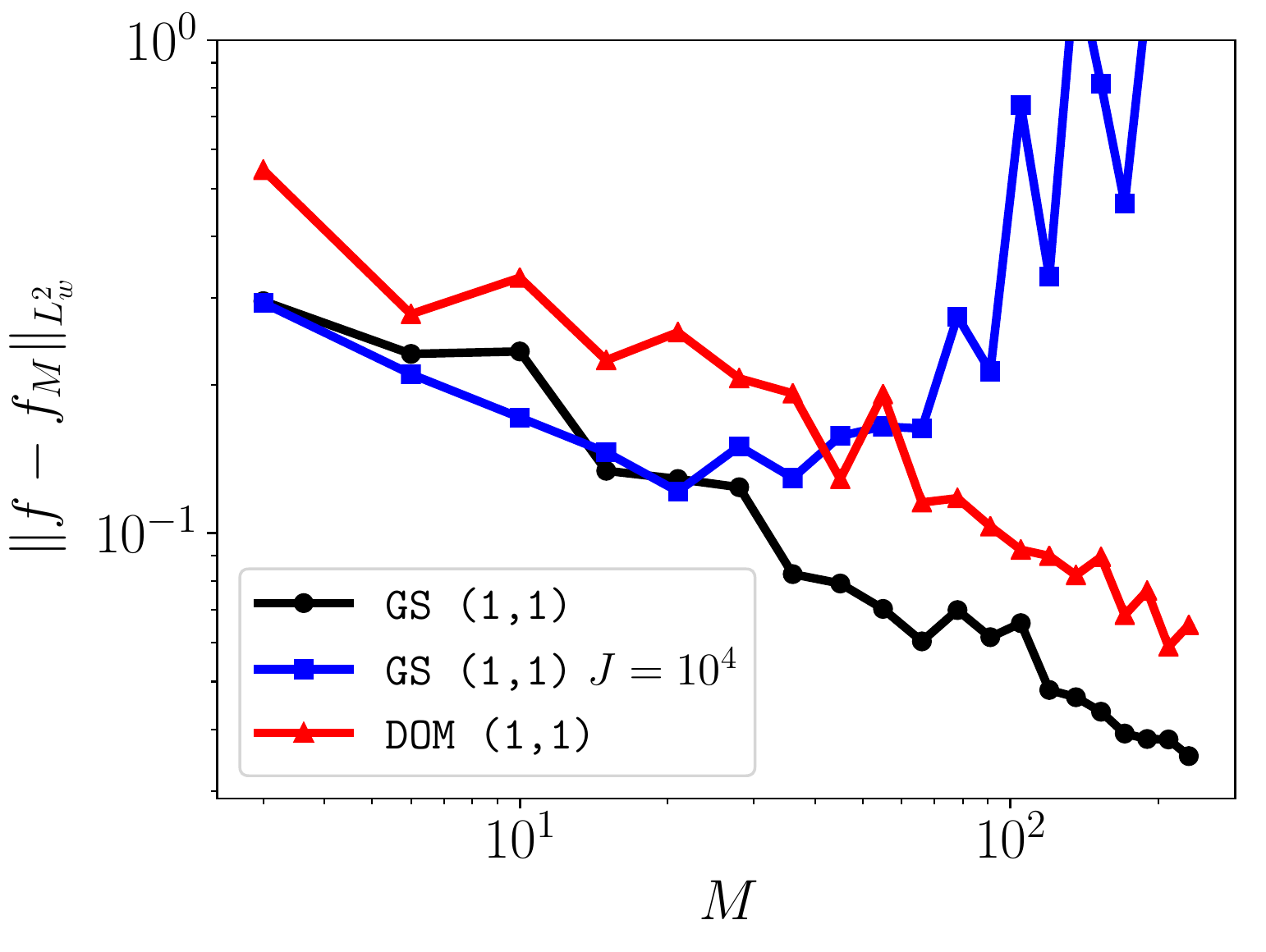}
\caption{(Left) Contour plot of the chemical reaction model response and the contours (lines) of the Banana density \eqref{eq:banana-density}. (Right)   \(L^2_\pbwt\)-errors in the polynomial interpolants of the chemical reaction model using the banana density. Each curve represents the median over 10 trials. $\ngs$ refers to the number of Monte Carlo samples used to compute the GS basis. When no $\ngs$ is given \eqref{eq:gs-tp-quadrature} was used. In the legend of the right plot we dropped the notation \texttt{($\alpha,\beta$)} because we used monomials for this example instead of Jacobi polynomials.}
\label{fig:banana-density}
\end{figure}

In the following we will assume that the rate parameters \(a=2(\rv_1+3)/3\in[0,4]\) and \(b=30(\rv_2+2)/7+5\in[5,35]\) in \eqref{eq:chemical-species} are random variables where the posterior distribution \eqref{eq:posterior} of \(\brv\) is
\begin{align}
\label{eq:banana-density}
  \pbwt(\rvd) &= C\exp(-(\frac{1}{10}\rvd_1^4 + \frac{1}{2}(2\rvd_2-\rvd_1^2)^2)), & {\rvd}&\in \dom = [-3,3]\times[-2,6],
\end{align}
where $C$ is a constant chosen to normalize $\pbwt$ so that it is a probability density. This density is called a ``banana density" and is a truncated non-linear transformation of a bivariate standard normal Gaussian distribution. 
The response surface of the mass fraction over the posterior density is depicted in the left of Figure \ref{fig:banana-density}. The response has a strong non-linearity which makes it ideal for testing high-order polynomial interpolation.

To approximate the response surface we use the \texttt{GS $(1,1)$} method. We compute the moments needed to orthogonalize the basis using Gauss-Legendre quadrature (including the weight in the integrand), and compare with the same approach using Monte Carlo quadrature. We use rejection sampling to draw independent samples from the banana density to avoid the effects of reduced effective sampling size that arise when using MCMC sampling. In the right of Figure \ref{fig:banana-density} we plot the \(L^2_\pbwt\)-errors for various total-degree polynomial interpolants (up to degree 20) of the chemical reaction model using the banana density. Similar to previous results, the GS basis obtained using exact moments achieves the smallest error for a fixed number of samples. Using samples from the density to estimate moments achieves comparable accuracy for small sample sizes until ill-conditioning resulting from the low accuracy of the moments takes effect.
\begin{figure}[htb]
\centering
\includegraphics[width=0.48\textwidth]{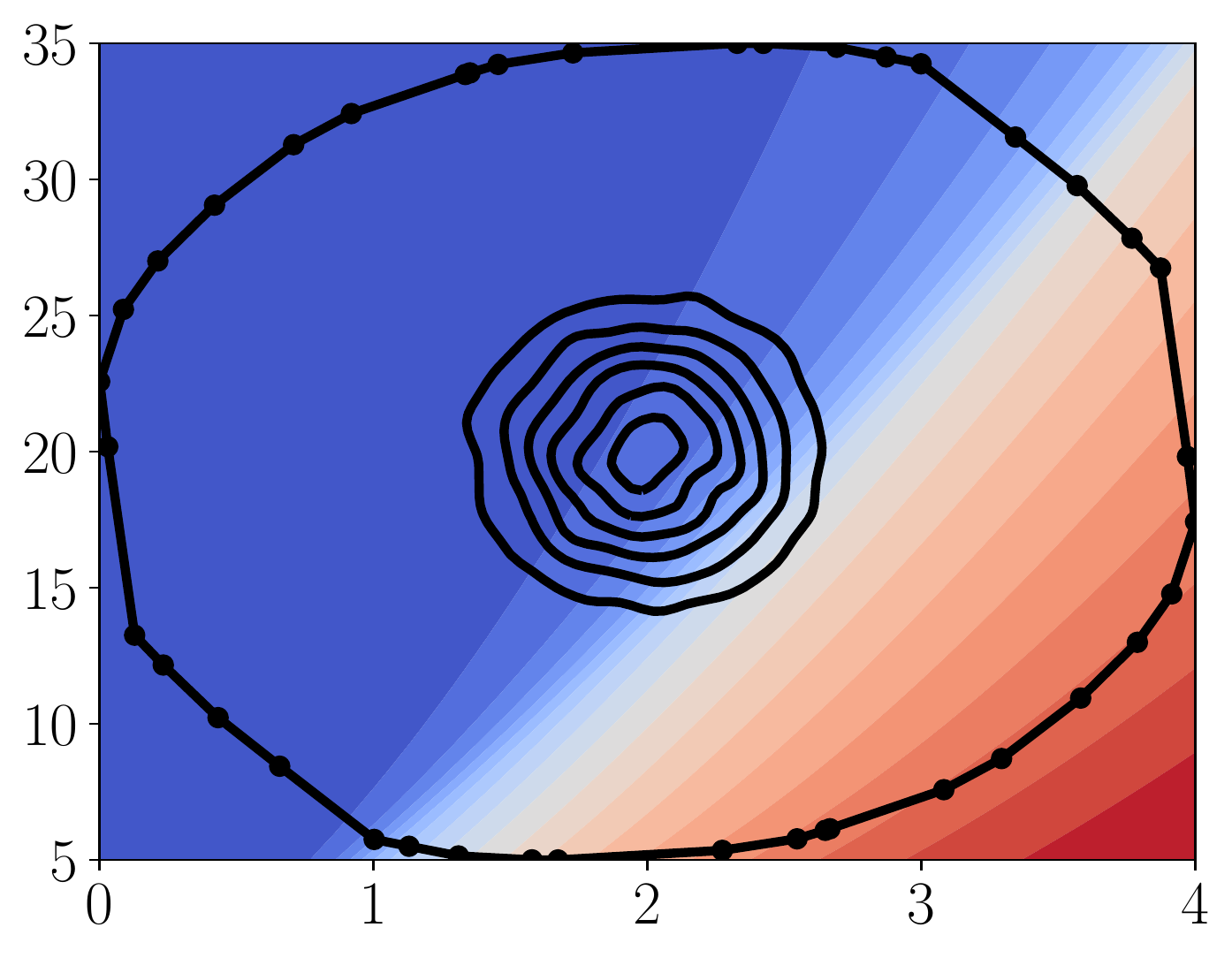}
\includegraphics[width=0.48\textwidth]{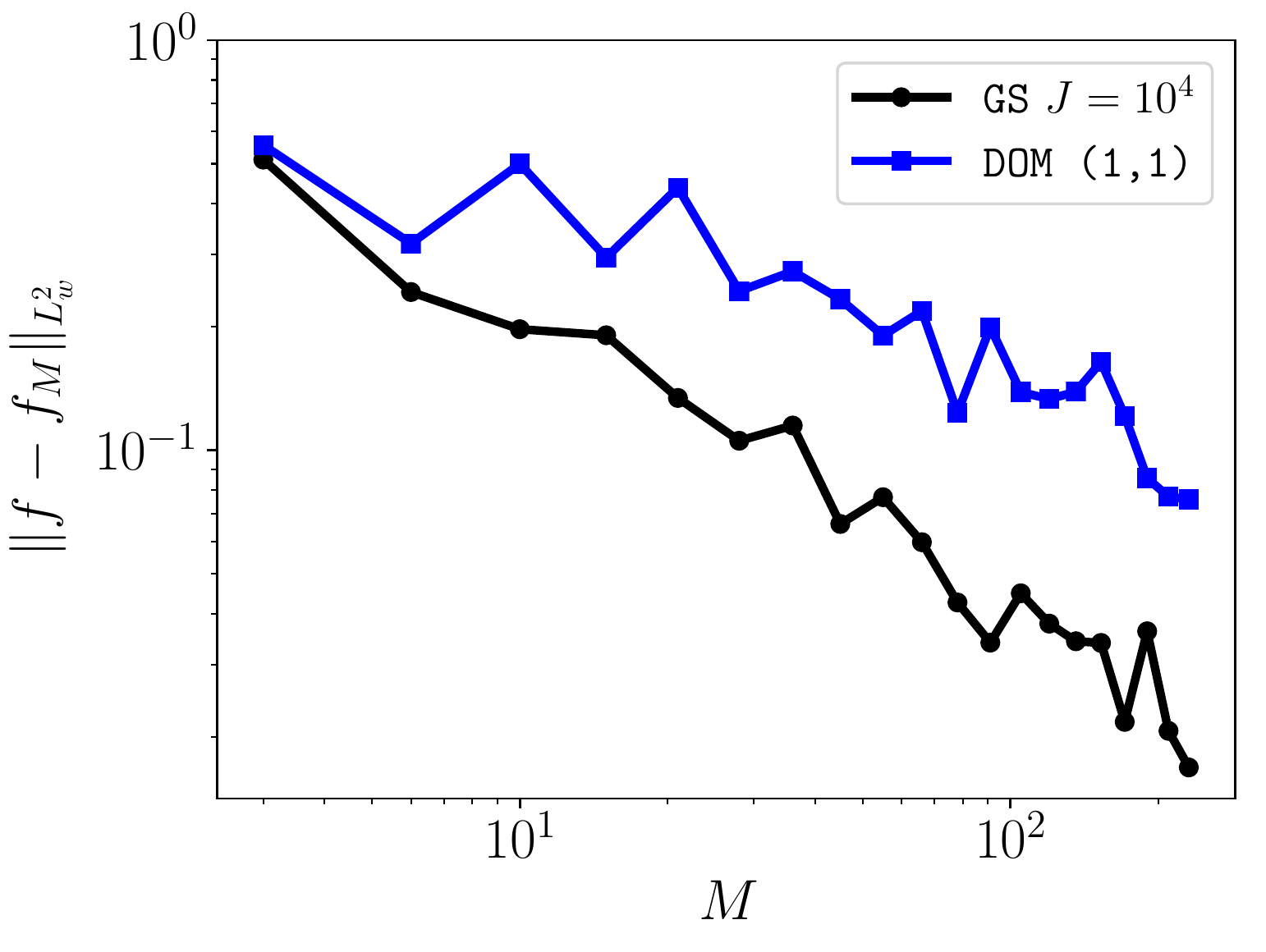}
\caption{
  (Left) Contour plot of the chemical reaction model response and the contours (lines) of the zonotope density. \(L^2_\pbwt\)-errors in the polynomial interpolants of the chemical reaction model using the zonotope density (right). Each curve represents the median over 10 trials. $\ngs$ refers to the number of Monte Carlo samples used to compute the GS basis. When no $\ngs$ is given \eqref{eq:gs-tp-quadrature} was used.}
\label{fig:chemical-reaction-convergence}
\end{figure}

We remark that previous attempts at building an approximation to increase the efficiency of Bayesian inference can be found in~\cite{Conrad_MPS_JASA_2016,Li_M_SISC_2014,Mattis_W_CMAME_2018}. Two of these methods~\cite{Conrad_MPS_JASA_2016,Mattis_W_CMAME_2018} use low-order localized surrogates to facilitate sampling in regions of high-probability. The use of localized surrogates results in small-rates of convergence. This is in contrast to the method proposed in this paper, which can exhibit spectral convergence rates. The method presented in~\cite{Li_M_SISC_2014} uses a sequence of global approximations that eventually concentrate on the posterior distribution. However unlike our proposed approach, the method in~\cite{Li_M_SISC_2014} does not use all model evaluations in the final approximation. Many samples are only used to construct the intermediate approximations.

\subsubsection{Dimension reduction}
\label{sec-5-5-2}
Using the same chemical species problem \eqref{eq:chemical-species}, we show here that our algorithm can be used to integrate high-dimensional ridge functions.
Ridge functions are multivariate functions that can be expressed as a function of a small number of linear combinations of the input variables variables \cite{Pinkus_book_2015}. We define a ridge function to be a function of the form 
\(f: \mathbb{R}^d\rightarrow\mathbb{R}\) that can be expressed as a function \(g\) of \(s < d\) rotated variables, 
\begin{align*}
  f(y) &= g(\rotation y), & \rotation &\in\mathbb{R}^{s\times d}.
\end{align*}
In applications it is common for \(d\) to be very large, but for \(f\) to be (an approximate) ridge function with \(s \ll d\). When approximating a ridge function one does not need to build a surrogate in \(\mathbb{R}^d\) but rather can focus on the more tractable problem of approximating in \(\mathbb{R}^s\). The difficulty then becomes generating an orthonormal basis with respect to the low-dimensional transformed probability density, which is typically unknown and challenging to compute. When the \(d\)-dimensional space is a hypercube, then the corresponding \(s\)-dimensional probability density is defined on a multivariate zonotope, i.e., a convex, centrally symmetric polytope that is the \(s\)-dimensional linear projection of a \(d\)-dimensional hypercube. The vertices of the zonotope are a subset of the vertices of the \(d\)-dimensional hypercube projected onto the \(s\)-dimensional space via the matrix \(\rotation \). 

In the following we will consider the approximation of a ridge function of independent uniform random variables with uniform density \(v(y)\). We set the ridge function to be a function of the mass fraction of the first species \(u_1\) of the competing species model \eqref{eq:chemical-species}. Specifically we define $\rvd=\rotation y \in \R^2$, where $y\in\mathbb{R}^{20}$ and $\rotation \in\mathbb{R}^{2\times 20}$ is a randomly generated matrix with orthogonal rows, and we define the ridge function to be \(f(y)=u_1(\rotation y)\).  The projection of the independent high-dimensional density \(v\) into the lower dimensional space induces a new dependent probability density $\pbwt$ on a zonotope $\Omega$. The resulting zonotope and probability density are shown in the left of Figure \ref{fig:chemical-reaction-convergence}. The zonotope density is obtained from a KDE constructed using \(10^4\) samples \(\brv=\rotation Y\), where \(Y\) were drawn randomly from the \(d=20\) uniform distribution.

To build an orthonormal basis on the zonotope we must select a tensor-product basis, a quadrature rule to compute the moments of that basis, and a procedure for generating candidate samples for the Leja sequence. Here we use a monomial basis and Sobol sequences to compute moments of the monomial basis with respect to \(\pbwt\). Specifically we generate a Sobol sequence of \(J=10^4\) \(d=20\)-samples \(\{y^{(j)}\}_{j=1}^J\) and then project these samples onto the zonotope to generate a quadrature rule \(\mathcal{\rvd}^J,w\), where \(\rvd^{(j)}=A y^{(j)}\) and \(w_q=1/J\). To generate candidate samples we use a randomized algorithm \cite{Stinson_GC_ARXIV_2016} to find the vertices of the zonotope and use rejection sampling to sample inside the convex hull of the zonotope vertices. 

In the right of Figure \ref{fig:chemical-reaction-convergence} we plot the \(L^2_\pbwt\)-errors in various polynomial total-degree interpolants (up to degree 20) of the chemical reaction model ridge function. We compare the GS procedure against a \texttt{DOM $(1,1)$} procedure. Again the GS basis achieves a significantly smaller error, for a fixed number of samples, than when using a tensor product basis. In this example the accuracy of the moments obtained by Sobol quadrature is sufficient to avoid noticeable effects of ill-conditioning for the ranges of degrees considered.
\rev{
\subsection{Partial differential equations with random input data}
In this section we apply our proposed methodology to build an approximation of a QoI obtained from a model of diffusion in porous media. We utilize this example to explore the effect of dimensionality on our proposed methodology.

Consider the following model of diffusion
\begin{align}
\label{eq:diffusion}
\frac{d}{dx}\left[k(x,\rvd)\frac{du}{dx}(x,\rvd)\right] = \cos(x,\rvd)& &
(x,\rvd)\in(0,1)\times \dom \\\nonumber
u(0,\rvd)=0 \quad  u(1,\rvd)=0
\end{align}
where the random diffusivity $k$ is a random field represented by the expansion
\begin{align}
\label{eq:diffusivityZ}
\log(k(x,\rvd)-0.5)=1+\rvd_1\left(\frac{\sqrt{\pi L}}{2}\right)^{1/2}+\sum_{k=2}^d \lambda_k\xi(x)\rvd_k,
\end{align}
where
\begin{align}
  \lambda_k=\left(\sqrt{\pi L}\right)^{1/2}\exp\left(-\frac{(\lfloor\frac{d}{2}\rfloor\pi L)^2}{8}\right) k>1,  & &  \xi(x)=
    \begin{cases}
      \sin\left(\frac{(\lfloor\frac{d}{2}\rfloor\pi x)}{L_p}\right) & k \text{ even}\\
      \cos\left(\frac{(\lfloor\frac{d}{2}\rfloor\pi x)}{L_p}\right) & k \text{ odd}
    \end{cases}
\end{align}
Given a correlation length $L_c$, which controls the variability of the random field, we set $L_p=\max(1,2L_c)$ and $L=\frac{L_c}{L_p}$.

The major challenge of polynomial approximation is the fast growth of the polynomial basis with the number of random variables. To mitigate this issue, we use an anisoptropic index set that consists of higher degree terms in important directions and lower-degree approximations in the less important directions. Specifically we use the following index sets presented in~\cite{Nobile_TW_SIAMNA_2008}:
\begin{align}\label{eq:anisotropic-indices}
  \Lambda_\alpha(l)=\bigcup_{\gamma\in\Gamma_\alpha(l)} \{\lambda\mid \lambda_k \le \gamma_k,k=1,\ldots d \},  & & \Gamma_\alpha(l)=\left\{ \gamma\in\mathbb{N}_{+}^d,\gamma_k\ge 1 \mid \sum_{k=1}^d (\gamma_k-1)\alpha_k\le l \alpha_\text{min} \right\},
\end{align}
where $\alpha_\text{min}=\min_k \alpha_k$ and
\begin{align}
  \alpha_k=
  \begin{cases}
    \frac{1}{2}\log\left(1+\sqrt{\frac{1}{24\sqrt{\pi}L}}\right) & k=1\\
    \frac{1}{2}\log\left(1+\sqrt{\frac{1}{48\sqrt{\pi}L}}\right)\exp\left(\frac{(\lfloor\frac{d}{2}\rfloor\pi L)^2}{8}\right) & k>1
  \end{cases}
\end{align}

To study the convergence of the approximation methods discussed in this paper, we consider problems of varying number of variables $d$ and investigate the convergence of the $L^2_\pbwt$ error as the number of points used to build each approximation increases. With this goal, we set $L_c=0.5$ and the joint density of $\rv$ to be the mixture of two independent tensor-product beta distributions\footnote{Note that, although the individual components of the mixture are independent, the mixture of the components is not.  Indeed, the distribution is bi-modal.}, i.e. $\pbwt(\rvd)=\frac{1}{2}B(\rvd,10,4)+B(\rvd,4,10)$. In Figure~\ref{fig:apc-diffusion-example}, we plot the error (computed using 1000 validation samples) as a function of the number of Leja samples $M$, which is dependent on $l$ in~\eqref{eq:anisotropic-indices}.

\begin{figure}[htb]
\centering
\includegraphics[width=\textwidth]{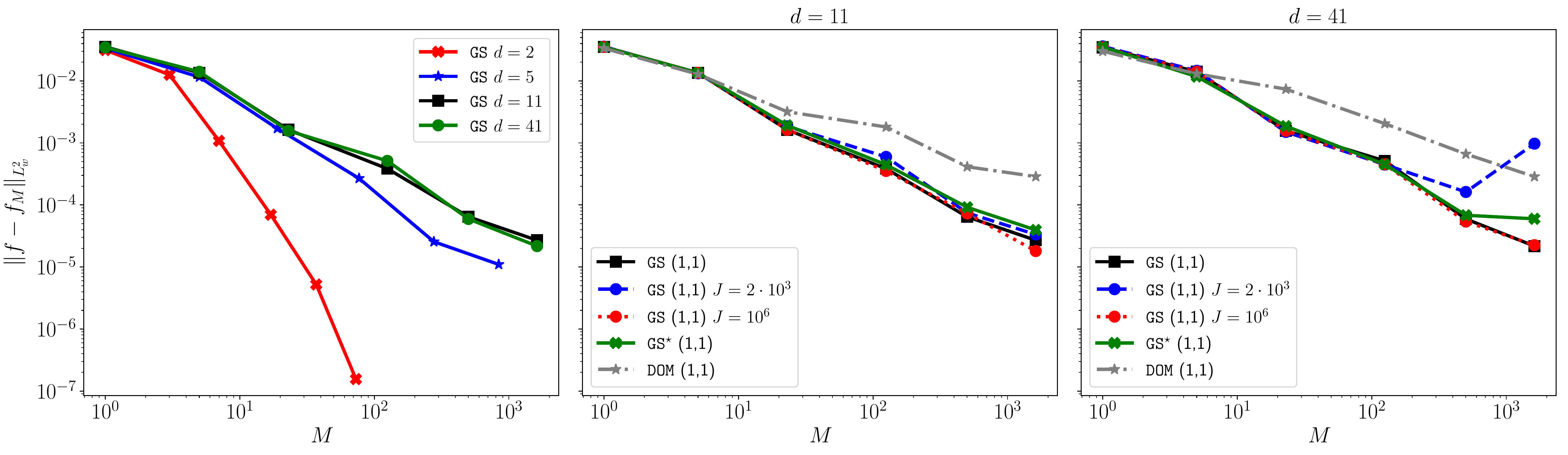}
\caption{\rev{A comparison of the $L^2_\pbwt$ error in the PCE approximations of the QoI $u(0.5,\rv)$ obtained using \eqref{eq:diffusion} and diffusivity \eqref{eq:diffusivityZ} with $L_c=0.5$. In the left plot, we display the error in $\texttt{GS}$ polynomials, as the number of samples and variables is increased. In the middle and right plots, we fix the dimensionality at $d=11$ and $d=41$, respectively, and compare errors in $\texttt{GS}$ polynomials orthogonalized with different quadrature rules with errors for $\texttt{DOM}$. $\ngs$ refers to the number of Monte Carlo samples used to compute the GS basis. When no $\ngs$ is shown, we use sparse grid quadrature to compute the inner products of the GSO procedure. All GSO Leja sequences were built using 10,000 candidate samples, except $\texttt{GS}^\star$ which only used 1,000.}}
\label{fig:apc-diffusion-example}
\end{figure}

From the left plot, we see that the rate of convergence of our $\texttt{GS}$ approach\footnote{We use sparse grid quadrature to produce an orthonormal basis. The index set~\eqref{eq:anisotropic-indices} can be exactly integrated by sparse grids for any tensor-product probability measure using weighted Leja sequences~\cite{Narayan_J_SISC_2014}. When $\pbwt$ is a mixture of independent measures, we can simply use a separate sparse grid quadrature rule for each mixture component and estimate the true mixture integral using a weighted sum of the individual quadrature rules.} does not degrade as the number random variables increases for $d>11$. This is consistent with the behavior observed for sparse grid approximations (with i.i.d uniform variables) in~\cite{Nobile_TW_SIAMNA_2008}. In the middle and left plots, we observe that the $\texttt{GS}$-based PCE can be significantly more accurate than the $\texttt{DOM}$-based PCE. The accuracy of the $\texttt{GS}$ PCE is degraded, relative to $\text{DOM}$, only when the accuracy of the quadrature rule is insufficient to accurately othonormalize the tensor product basis used in~\eqref{eq:gram-schmidt-alg} or when insufficient candidates are used to build the Leja sequence. There is a mild dependence on dimension for both the number of Monte Carlo samples needed to produce a well-conditioned GSO basis and the number of candidate samples used to construct a conditioned Leja sequence. These two effects are evident by the deterioration of the error in the $\texttt{GS}$ (1,1) $J=2\cdot 10^3$ and $\texttt{GS}^\star$ approximations, respectively, when the number of variables is increased from 11 to 41.}

\section{Conclusions}
\label{sec-6}
\rev{In this paper, we have presented an algorithm for building polynomial interpolants of functions of dependent variables. Most existing literature focuses on strategies for building polynomial approximations under the assumption of independence. Our work provides two major contributions to the existing literature, namely an investigation of the disadvantages of existing polynomial interpolation methods for approximation with dependent variables and an algorithm for generating interpolants that minimizes the ill-conditioning of the interpolation matrix.

We demonstrate through extensive numerical examples that our results are almost always significantly more accurate than existing approaches. The improved performance is obtained by balancing the need to sample in high-dimensions with minimizing the ill-conditioning of the polynomial interpolation matrix. The efficacy of our approach is only limited by the ability to produce multivariate polynomials which are orthonormal to the joint distribution of the random variables. We use Gram-Schmidt orthogonalization to generate such polynomials. The stability of this procedure is dependent on the dependence structure between the random variables and the accuracy of the quadrature rules used to compute the inner-products in the Gram-Schmidt procedure.

The ill-conditioning of Gram-Schmidt introduced by using approximate numerical integration can be minimized by using highly accurate quadrature rules, provided an explicit form for the joint density is available. However, in a number of practical situations, one does not know the joint density explicitly, but rather only has samples drawn from its distribution. In these settings, for example when Bayesian inference is used to condition prior estimates of uncertainty on data, we can still compute highly accurate interpolants, but the accuracy of our approach is dependent on the number of available distribution samples. In future work, we will investigate alternative strategies for producing orthonormal multivariate polynomials that do not suffer from the ill-conditioning associated with the Gram-Schmidt procedure.

Our approach is extremely flexible and can exploit various types of structure in the function being approximated. We demonstrated this by leveraging dimension reduction and {\it a priori} estimates of anisotropy to build interpolants that can be used for high-dimensional approximation. In future work, we plan to develop adaptive strategies for determining the best polynomial index set that minimizes cost while maximizing accuracy, when anisotropy must instead by discovered.}

\bibliographystyle{plain}

\bibliography{references}

\section{Acknowledgements}
\label{sec-7}
J. Jakeman and M. Eldred were supported by DARPA EQUiPS and DOE SCIDAC. Sandia National Laboratories is a multi-mission laboratory managed and operated by National Technology and Engineering Solutions of Sandia, LLC., a wholly owned subsidiary of Honeywell International, Inc., for the U.S. Department of Energy's National Nuclear Security Administration under contract DE-NA-0003525. The views expressed in the article do not necessarily represent the views of the U.S. Department of Energy or the United States Government.

A. Narayan is partially supported by NSF DMS-1720416, AFOSR FA9550-15-1-0467, and DARPA EQUiPS N660011524053.

F. Franzelin and D. Pfl\"{u}ger were supported by the German Reseach Foundation (DFG) within the Cluster of Excellence in Simulation Technology (EXC 310).

\section{Appendix}
\label{sec-8}

\subsection{Generating samples using the Gaussian copula}
\label{sec-8-2}
The following algorithm can be used to generate dependent samples, from a Gaussian copula with correlation matrix \(\Rvy\), with marginals,  with marginal probability densities given by \(f_i(z_i)\). This algorithm can be used to generate random variables which satisfy the assumptions of the Nataf transformation~\ref{sec-4-2}.
\begin{enumerate}
\item Generate a sample \(u=(u_1,\ldots, u_d)\) from the multivariate Gaussian distribution \(N(0,I)\) with zero mean and unit variance.
\item Compute the Cholesky factorization \(L\) of the correlation matrix such that \(\Rvy=LL^T\).
\item Compute a correlated standard normal sample \(v=Lu\).
\item Modify the sample \(v\) to have uniform marginals in every dimension. I.e. compute \(\hat{v}_i=\Phi(v_i)\), where \(\Phi\) is the CDF of the standard normal distribution.
\item Generate a sample \(z=(z_1,\ldots,z_d)\) from the desired distribution. I.e. compute \(z_i=F_i(\hat{v}_i)\), where \(F_i\) is the CDF of the desired marginal density.
\end{enumerate}
\end{document}